\documentclass[onefignum,onetabnum]{siamart190516}

\usepackage{lipsum}
\usepackage{amsfonts}
\usepackage{graphicx}
\usepackage{epstopdf}
\usepackage{algorithmic}
\Crefname{ALC@unique}{Line}{Lines}
\ifpdf
  \DeclareGraphicsExtensions{.eps,.pdf,.png,.jpg}
\else
  \DeclareGraphicsExtensions{.eps}
\fi

\newtheorem{thm}{Theorem}
\newtheorem{prop}[thm]{Proposition}

\newtheorem{rem}{Remark}
\newtheorem{ex}{Example}
\newtheorem{defn}[thm]{Definition}

\newcommand{\dd}{{\,\mathrm d}}
\newcommand{\si}{\sigma}

\renewcommand{\th}{\theta}
\newcommand{\la}{\lambda}
\newcommand{\ga}{\gamma}

\newcommand{\eps}{\varepsilon}

\renewcommand{\phi}{\varphi}

\newcommand{\scr}[1]{{\mathcal #1}}
\newcommand{\argmin}{\operatornamewithlimits{argmin}}

\newcommand{\PP}{\mathbb{P}}
\newcommand{\ind}{\mathbf{1}}
\newcommand{\RR}{\mathbb{R}}

\newcommand{\Bm}{\begin{bmatrix}}
\newcommand{\Em}{\end{bmatrix}}
\newcommand{\bs}[1]{{\boldsymbol #1}}
\newcommand{\T}{{\prime}}

\newcommand{\Md}{M^\dagger}
\newcommand{\GP}{\ensuremath{\scr{G}\scr{P}}}

\newcommand{\innerproduct}[1]{\left<#1\right>} 
\newcommand{\Diff}{\operatornamewithlimits{Diff}}
\newcommand{\Id}{\operatornamewithlimits{Id}}

\newsiamremark{remark}{Remark}
\newsiamremark{hypothesis}{Hypothesis}
\crefname{hypothesis}{Hypothesis}{Hypotheses}
\newsiamthm{claim}{Claim}

\headers{Diffusion bridges for stochastic Hamiltonian systems}{A. Arnaudon, F. v d Meulen, M. Schauer, and S. Sommer}

\title{Diffusion bridges for stochastic Hamiltonian systems and shape evolutions\thanks{Submitted to the editors DATE.
}}

\author{Alexis Arnaudon\thanks{Imperial College London, UK}
\and Frank van der Meulen\thanks{Delft University of Technology, NL.}
\and Moritz Schauer\thanks{Chalmers University of Technology, SE.}
\and Stefan Sommer\thanks{University of Copenhagen, DK.}
}

\usepackage{amsopn}

\ifpdf
\hypersetup{
  pdftitle={Diffusion bridges for stochastic Hamiltonian systems and shape evolutions},
  pdfauthor={A. Arnaudon, F. v d Meulen, M. Schauer, and S. Sommer}
}
\fi

\begin{document}

\maketitle

\begin{abstract}
Stochastically evolving geometric systems are studied in shape analysis and computational anatomy for modelling random evolutions of human organ shapes. The notion of geodesic paths between shapes is central to shape analysis and has a natural generalisation as diffusion bridges in a stochastic setting. Simulation of such bridges is key to solve inference and registration problems in shape analysis. We demonstrate how to apply state-of-the-art diffusion bridge simulation methods to recently introduced stochastic shape deformation models thereby substantially expanding the applicability of such models. We exemplify these methods by estimating template shapes from observed shape configurations while simultaneously learning model parameters.

\end{abstract}



\section{Introduction}
We develop a general scheme for simulation of conditioned diffusion processes, diffusion bridges, for finite dimensional landmarks manifolds appearing in fields including shape analysis and fluid dynamics. 
We demonstrate several important properties of the scheme, including that the approach applies to more general classes of stochastic models than presently covered in the literature, that it effectively  handles  challenging aspects of shape spaces, nonlinearity and partial observability, and that it generalizes the inexact matching scheme that is commonly used in applications to the stochastic setting. 
Importantly, the simulation scheme allows for statistical inference of properties of the evolution of the shapes. We exemplify this by estimating template shapes from observed shape configurations.

\subsection{Background}

The geometry of landmark spaces and the landmark matching problem, exemplified by matching of finite sets of anatomical markers on medical images, have been extensively studied in shape analysis.  Examples include the seminal work of Kendall~\cite{kendall_shape_1984} on landmark configurations modulo scaling and rotation, and the extensive interest of landmark dynamics for diffeomorphic shape models~\cite{joshi_landmark_2000,younes_shapes_2010}. The landmark matching problem has a geometric formulation in the language of geometric mechanics where diffeomorphisms of the underlying domain act to move the landmarks, and optimal trajectories between observed landmark configurations satisfy geodesic equations for a right-invariant metric on the diffeomorphism group (\cite{younes_evolutions_2009}). These equations are also called Euler-Poincar\'e equations~\cite{holm_euler-poincare_1998}. 
Landmark dynamics are intrinsically linked to soliton dynamics in fluid dynamics~\cite{holm_soliton_2004}, e.g.\  soliton solutions of the Camassa-Holm equations~\cite{camassa_integrable_1993}, or more generally singular solutions of the Euler-Poincar\'e equation~\cite{holm2005momentum}. Stochastic models of landmark evolutions are therefore of interest in both shape analysis and fluid dynamics. 

Several  recent models include stochastic effects in landmark analysis~\cite{markussen_large_2007,trouve_shape_2012,vialard_extension_2013,Marsland-Shardlow,arnaudon_stochastic_2017,arnaudon_geometric_2019}. The interest appears for several reasons: 
\begin{enumerate}
\item
 If $v$ denotes the observed landmark configuration, the presence of noise in the observation implies that landmark analysis is classically performed in the inexact setting where the requirement that the flow $q_t$ satisfy $q_T=v$ at observation time $t=T$ is essentially relaxed to $q_T=v+\epsilon$ for some perturbation $\epsilon$. This thus amounts to adding noise at the endpoint $q_T$ of the flow $q_t$. If the perturbation is not linked to observation noise but instead comes from intrinsic stochastic behaviour of the shape or soliton evolution, it is more natural to make the noise intrinsically time-continuous, i.e.\ added to $q_t$ for each $t$ on an infinitesimal level.
 \item  The flow $q_t$ can have multi-scale behaviour where the coarse scale evolution is modeled by a deterministic flow while the fine scale behaviour, which in a fluid system is closer to turbulence, can be summarised with stochastic terms. 
 \item  In medical applications, it can be hypothesised that the shape evolution is not purely deterministic, and that the stochastic deformation happen continuously in time. 
 \item  It is generally hard to construct natural families of probability distributions on geometric spaces without linearising the spaces around a center point. Solutions to stochastic differential equations however provide such probability models \cite{sommer_probabilistic_2020}.
\end{enumerate}

In this paper, we focus on models defined by stochastic differential equations. In  all of the above enumerated cases, it is common that observations are  discrete in  time. Constructing paths that are conditioned on such partial in time observations is called  {\it bridge simulation} and the paths {\it bridges}. Bridge simulation is essential for likelihood based inference as it provides a stochastic method for approximating the intractable likelihood. If possible, it enables methods such as sequential Monte Carlo, Markov Chain Monte Carlo, stochastic Expectation-Maximisation, or stochastic gradient descent. This in turn allows for estimation of parameters in the process' dynamics and recovery of latent states at observation times. Moreover, in applications like landmark matching, there is intrinsic interest in bridge simulation as bridges appear as stochastic perturbations of landmark geodesics or as critical paths for a stochastically perturbed energy. Bridges therefore present a stochastic generalisation of the geodesic paths usually considered when matching sets of landmarks.

Bridge simulation is generally a challenging problem. It has been addressed in the Euclidean setting in a series of works including~\cite{DelyonHu,beskos_exact_2006,PapaspiliopoulosRoberts,bladt_simulation_2016,schauer_guided_2017}.
For shape manifolds, previous papers have centred on simulation in landmark spaces. Manifolds of finite numbers of landmarks are generally of finite dimensionality in contrast to spaces of continuous shapes which are in essence infinite dimensional. While bridge sampling on infinite dimensional shape spaces remains an open question, bridge sampling for landmarks spaces has been considered recently~\cite{arnaudon_stochastic_2017,arnaudon_geometric_2019,sommer_bridge_2017} using the approach of~\cite{DelyonHu,marchand_conditioning_2011}. 
The underlying diffusion arises from stochastic Euler-Poincar\'e equations~\cite{holm_variational_2015} and, in a different setting, Brownian motion on the landmark space~\cite{sommer_bridge_2017}. The difficulty in landmark bridge sampling arises from the high-dimensionality, the curvature and hence nonlinearity of the diffusion process, the degenerate diffusion processes in the stochastic Euler-Poincar\'e case (vanishing eigenvalues in the diffusion matrix), and from observations being only in the position variable of the phase-space flow.
The methods from~\cite{DelyonHu} and~\cite{arnaudon_stochastic_2017,arnaudon_geometric_2019} apply only for certain classes of stochastic evolution, excluding other important stochastic landmark approaches such as proposed in~\cite{trouve_shape_2012,Marsland-Shardlow}, due to the degeneracy of the flow. As stated in~\cite{Marsland-Shardlow}: \emph{``Unfortunately, none of the known methods for diffusion bridges works with (2.1) to give computationally convenient algorithms. Without an efficient method for sampling the diffusion bridge, it is hard to formulate an MCMC method with good acceptance rates. Consequently, the generalised Langevin prior distribution is difficult to use in Bayesian statistics, and we now turn to simpler prior distributions, which arise by approximating the Langevin equation.''} In the present paper, we show that resorting to such simpler prior distributions is unnecessary. 

\subsection{Relation to image registration}
Landmark matching in the Large Deformation Diffeomorphic Metric Mapping (LDDMM, \cite{younes_shapes_2010}) setting is closely related to image registration, both because landmarks often arise from annotations of images, and because the LDDMM model allows to generate deformations of the full image domain from matching of only a small number of landmarks. As such, landmark matching can partly be considered a reduced problem compared to matching of full images. In this paper, we solve a different, but related, problem compared to image registration: we aim to make possible statistical inference of properties of the evolution of shapes, specifically in the form of landmark configurations. The aim is thus not to find a good diffeomorphic matching between two images or shapes, but to infer statistical properties of stochastic transformations between shapes. Nevertheless, some of the stochastic models considered in the paper, particularly the Eulerian model, allows to generate full domain deformations as is the case for LDDMM landmark matching (this is called \emph{lifting} of the trajectory from the landmark space to the space of diffeomorphisms). One can therefore envision making bridge sampling with data sizes - full images or large numbers of landmarks - comparable to what is currently used in image registrations algorithm. While this will be extremely interesting, it is also very challenging both from a computational perspective and theoretically due to the high or even infinite dimensionality of the resulting shape space. In the present paper, we aim to present one of the first applications of bridge sampling in shape analysis that nevertheless applies bridge sampling techniques to stochastic processes of much higher dimensionality than what is regularly seen in other areas of statistics. We will pursue extension to deformations of the full domain in future work.

\subsection{Contribution}
We extend the Euclidean simulation approach of~\cite{schauer_guided_2017,partialbridges} to a general method for simulating stochastic landmark equations. We will pursue this under the stochastic models for landmark and soliton dynamics introduced in~\cite{trouve_shape_2012}, \cite{Marsland-Shardlow} and~\cite{arnaudon_stochastic_2017}. 
The presented method is the first in the literature to allow bridge sampling for the former schemes, and it improves upon previous sampling schemes by better incorporating nonlinearity in the drift and diffusion coefficients of the SDEs. For models with intrinsic noise, it is the first  approach that provides full uncertainty quantification via samples of the posterior. The approach generalises the notion of inexact shape matching as in~\cite{younes_shapes_2010}. We demonstrate how the sampling approach can be used to estimate the starting conditions of the diffusion and parameters of the diffusion dynamics. Particularly, we apply it for the problem of estimating a template shape given observed shape configurations. In this case, we show how using the Riemannian Manifold Metropolis Adjusted Langevin Algorithm (RMMALA) is crucial to obtain efficient updates of the estimated template.

\subsection{Outline}

The paper starts in Section~\ref{sec:review} with a review of finite dimensional shape analysis and its stochastic extensions. In Section~\ref{sec:bridge} we explain how guided proposals can be used for bridge simulation. Their specific implementation for the two stochastic landmarks models considered in this work is specified in Section~\ref{sec:guided_ms_ahs}.  Markov Chain Monte Carlo  sampling schemes are detailed in Section~\ref{sec:mcmc} and   experimental evaluation on synthetic and empirical datasets are given in Section~\ref{sec:numexamples}. We end with a discussion section with suggestions for future work. 

\section{Finite Dimensional Shape Analysis}\label{sec:review} 

To provide the necessary background for the stochastic landmark models, we here give a short review of landmark shape spaces with geometry inherited from a right-invariant metric on the diffeomorphism group.
Landmark shape spaces are finite dimensional in contrast to spaces of continuous curves and surfaces that are inherently infinite dimensional. On a formal level, the geometric setting as outlined below is however equivalent. This makes the landmark case an example of how stochastics and bridge simulation can be used in a geometric setting while we can postpone the intricacies of infinite dimensionality in other shape spaces to later works.

Shape analysis as pursued in the Large Deformation Diffeomorphic Metric Mapping (LDDMM, \cite{younes_shapes_2010}) framework starts with actions of the diffeomorphism group on shapes spaces. For landmarks, let $q=(q_1,\ldots,q_n)$ be a configuration of $n$ distinct landmarks $q_i\in\Omega$ in a domain $\Omega\subset\mathbb R^d$. Let $\phi$ be an element of the set $\Diff(\Omega)$ of diffeomorphism on $\Omega$, smooth invertible mappings with smooth inverses. Then $\phi$ acts on $q$ by composition $\phi.q=(\phi(q_1),\ldots,\phi(q_n))$.
For fixed $q$, the map $\pi\colon\Diff(\Omega)\to\mathcal M$, $\pi(\phi)=\phi.q$ is denoted the action map.

The landmark space $\mathcal M=\{(q_1,\ldots,q_n)\in\Omega^n \mid q_i\not=q_j,\, i\not=j\}$ can be given the structure of a manifold by letting it inherit the differentiable structure from its embedding as an open subset of $\RR^{nd}$ ($\RR^{nd}$ except for the subset of points where landmark pairs coincide). It can furthermore be equipped with a Riemannian metric which in turn defines the length and energy of path in the landmark space, and from this the matching energy \eqref{eq:landmark_energy} used below is defined. This happens with the following geometric structure. 

{\bf 1)} Let $\mathcal X(\Omega)$ be the space of vector fields on $\Omega$. We equip subsets $V$ of  $\mathcal X(\Omega)$ with an inner product using a reproducing kernel Hilbert space (RKHS) structure: For $q_1, \dots, q_n \in\Omega$ let $V$ be the completion of the set of vector fields on the form $\sum_{i=1}^nK(\cdot,q_i)a_i$, $a_i\in\RR^d$. Here $K\colon \Omega\times\Omega\to\RR^{d\times d}$ is a matrix-valued map denoted a kernel.  $K$ is required to be symmetric and positive definite in the sense that 
\begin{equation*}
  \sum_{i,j=1}^n a_i^TK(q_i,q_j)b_j > 0
  \quad
  \forall a_i,b_j\in\RR^d,\ q_i,q_j\in\Omega,\ n\in\mathbb N
  \, .
\end{equation*}
In this case, $V$ receives an inner product by completing the inner product
\begin{equation}
  \innerproduct{K(\cdot,q_i)a_i,K(\cdot,q_j)b_j}_V
  =
  \sum_{i,j=1}^n a_i^TK(q_i,q_j)b_j
  \, ,
  \label{eq:Vinnerprod}
\end{equation}
to all of $V$ by continuity. This makes $V$ a Hilbert space.

{\bf 2)} Since $\mathcal{M}\subset\RR^{nd}$, tangent vectors $v\in T\mathcal{M}$ can be represented as vectors in $\RR^{nd}$. Let $q\in\mathcal{M}$ be a landmark configuration. The LDDMM Riemannian metric is then the inner product
\begin{equation}
  \innerproduct{v,w}_q
  =
  \sum_{i,j=1}^n v_i^TK(q_i,q_j)^{-1}w_j\, 
  \label{eq:metric}
\end{equation}
between vectors $v,w\in T_q\mathcal{M}$. Landmark dynamics are often described in terms of momenta, covectors in the dual bundle $T^*\mathcal{M}$. The corresponding inner product between momentum vectors $a,b\in T_q^*\mathcal{M}$ (the cometric) is then
\begin{equation}
  \innerproduct{a,b}_q
  =
  \sum_{i,j=1}^n a_i^TK(q_i,q_j)b_j
  \ ,
  \label{eq:cometric}
\end{equation}
i.e., the cometric coincides with the inner product \eqref{eq:Vinnerprod} on $V$. 

{\bf 3)} The Riemannian structure on $\mathcal{M}$ described here has its geometric origin in the following facts: The tangent space $T_{\Id}\Diff(\Omega)$ of the Lie group $\Diff(\Omega)$ is diffeomorphic to $\mathcal X(\Omega)$, and the Riemannian metric \eqref{eq:metric} is the Riemannian metric on $\mathcal{M}$ that makes the action map $\pi$ above a Riemannian submersion with respect to the right-invariant (invariant to the composition of diffeomorphisms on the right) Riemannian metric on $\Diff(\Omega)$ that the inner product \eqref{eq:Vinnerprod} on $V\subset\mathcal X(\Omega)$ specifies. In the current context, the main implication of this geometric argument is that the metric \eqref{eq:metric} and cometric \eqref{eq:cometric} have their definitions rooted on a well-defined geometric structure.

The Riemannian metric defines the energy of a path $\bs{q}=(q_t,\, t\in [0,T])$ of landmark configurations by
\begin{equation*}
  E(\bs q)=\int_0^T\innerproduct{\dot{q}_s,\dot{q}_s}_q \dd s
  \, .
\end{equation*}
This allows matching of landmark configurations $q$ and $v$ by searching for a minimal energy path $\bf q$, $q_t\in\mathcal M$ such that $q_0=q$ and $q_T=v$,
\begin{equation}
  \argmin_{\bs q,q_0=q,q_T=v}E(\bs q)
  \, ,
  \label{eq:landmark_energy}
\end{equation}
as pursued below.
Paths realising this energy are geodesics on $\mathcal{M}$. 

\subsection{Inverse problem of landmark matching}

The interest in landmark matching arises from the case where a fixed sets of meaningful landmarks, for example anatomical markers in medical images, can be identified. Shapes can then be analyzed by identifying differences between configurations of markers, in the medical case either between subjects or intra-subject when consecutive images are acquired over time. Interestingly, the LDDMM framework allows matching of landmarks as well as curves, surfaces and even images on $\Omega$ via the same framework with a metric on $\Diff(\Omega)$ descending to the particular shape spaces. 
However, we will restrict to the landmark case here for two primary reasons: The landmark manifold $\mathcal M$ is finite dimensional, and it has a natural embedding in $\RR^{nd}$ allowing us to write landmark dynamics in Euclidean coordinates.

Assume two landmark configurations $q$ and $v$ are given. A matching consists in solving the variational problem 
\begin{equation}
  \argmin_{\bs q}E(\bs q)
 \qquad 
  \text{such that}\qquad 
  q_0=q,\, q_1=v\, ,
  \label{eq:exact_match}
\end{equation}
with the energy $E$ from \eqref{eq:landmark_energy}. This is known as exact matching and results in a $t$-dependent transformation of the configuration $q$ to the configuration $v$. Due to the uncertainty often present in observations of landmarks, e.g. with manual annotations on medical images, inexact matching is instead often used in practice. The variational problem is now
\begin{align}
  \argmin_{\bs q}E(\bs q)+S(q_1,v)\qquad 
  \text{such that}\qquad 
  q_0=q\, ,
  \label{eq:inexact_match}
\end{align}
with the added term $S$ measuring the dissimilarity between the two landmark configurations $q_1$ and $v$. This term is often taken to be of the form $S(q_1,v)=n^{-1}\|q_1-v\|^2$ using the norm from the embedding of $\mathcal M$ in $\RR^{nd}$. Intuitively, we can think of Gaussian noise on the observations in $\RR^{nd}$ and $S$ being proportional to a log-likelihood.

Note that the $t$-dependence of solutions to \eqref{eq:exact_match} and \eqref{eq:inexact_match} allows us to think of $q_1$ as an observation of landmark trajectories $q_t$ at time $t=1$. This will in particular relate to the conditioning of stochastic processes later, where we assume a fixed observation time $T>0$. In the matching case, the notation implies $T=1$ without loss of generality.

\subsection{Deterministic landmark dynamics}

The energy \eqref{eq:landmark_energy} has a Hamiltonian formulation which involves the kernel $K$, where the Hamiltonian is
\begin{equation}
  H(q,p)
  =\frac12\innerproduct{p,p}_q
  =
  \frac12\sum_{i,j=1}^np_i^T K(q_i,q_j)p_j
  \, .
  \label{eq:Hamiltonian}
\end{equation}
Let $x \in \RR^{2dn}$ be the vector obtained by concatenating the tuples $\{(q_i, p_i)\}_{i=1}^n$. That is, $x$ represents an element of the phase space $T^*\mathcal M$ of the landmark space.
Let $K$ be a kernel as discussed above, and we assume it is scalar, i.e. of the form $K(y)=k(y)\Id_d$ for a real-valued function $k$. 
Hamilton's equations of motion are then explicitly given by
\begin{align}\label{hamiltonian}
    \begin{split}
      \frac{\dd}{\dd t} q_i  &=  \frac{\partial H}{\partial p_i} = \sum_{j=1}^n  p_j k(q_i-q_j)\, , \\
      \frac{\dd}{\dd t}  p_i  &= -\frac{\partial H}{\partial q_i}  =  -\sum_{j=1}^n \langle p_i, p_j \rangle  \nabla k(q_i-q_j)\, ,
    \end{split}
\end{align}
where $p$ is the momentum of the flow and  $\nabla k(y)$  the gradient of $k$.

Solutions $\mathbf q$ to \eqref{eq:exact_match} and \eqref{eq:inexact_match} appear as the $q$-variable of solutions $(\mathbf q,\mathbf p)$ to \eqref{hamiltonian}. Therefore, the search space for optimal solutions of the optimization problems can be reduced to solutions of Hamilton's equations, a reduction from an infinite to a finite dimensional space.

\subsection{Stochastic landmark dynamics}
As already indicated, there are several reasons to introduce stochastic evolutions.
We review here two different models of stochastic evolutions on landmarks which we will consider later. These are formulated in terms of stochastic differential equations. A good introduction to such equations is given in \cite{klebaner2012introduction}, whereas more advanced treatments are for example \cite{oksendal2003stochastic}, \cite{Karatzas-Shreve} and \cite{RogersWilliams2}. 

\subsubsection{Stochastic forcing: Lagrangian  and Langevin model}

The first stochastic perturbation of the dynamics landmarks described above has been introduced in~\cite{trouve_shape_2012,vialard_extension_2013} and consist in an additive Lagrangian noise of the form
\begin{align}
    \dd  q_i^\alpha &= \frac{\partial H}{\partial p_i^\alpha} \dd t \qquad \mathrm{and} \qquad \dd p_i^\alpha = -\frac{\partial H}{\partial q_i^\alpha} \dd t + \gamma_i \dd W_t^i\, , 
    \label{TV-model}
\end{align}
where $\gamma_i\in \RR$ is a scalar noise amplitude given for each landmark. We will refer to this system as the {\it Lagrangian model}. 

More recently, \cite{Marsland-Shardlow} generalised the Lagrangian model by adding a dissipative term to the model to obtain what we will call the Langevin model
\begin{equation}
	\dd p_i^\alpha =-\lambda \frac{\partial H}{\partial p_i^\alpha} \dd t   -\frac{\partial H}{\partial q_i^\alpha} \dd t + \gamma_i \dd W_t^i\, ,
	\label{TV}
\end{equation}
where $\la>0$ is a damping coefficient which implies existence of the Gibbs invariant measure for this SDE. This model is in the following referred to as the {\it Langevin model}.

\subsubsection{Transport noise: Eulerian model}

In addition to the Lagrangian noise of the previous model (where each Wiener process is associated to a landmark), \cite{arnaudon_geometric_2019} introduced an Eulerian noise, where the noise fields are functions of the domain $\Omega$. 
In general, any set of functions could be used, but in practice, we will fix a family of $J$ noise fields $\sigma_1, \ldots, \sigma_J$, which are kernel functions centred at locations $\{\delta_\ell\}$ of the form 
\begin{equation}\label{eq:noisefield}
    \si_\ell^\alpha(q) = \ga_\alpha\bar{k}_\tau(q-\delta_\ell)\, , 
\end{equation}
where $\bar{k}_\tau$ is a kernel with length-scale $\tau$ and $\gamma \in \RR^d$ is the vector of noise amplitudes. Note that it is possible to choose both $\gamma$ and $\tau$ dependent on the location $\delta_\ell$, though we will not use this flexibility here. 

The stochastic dynamics are then obtained as in the deterministic case, to be a Hamiltonian system but in both the drift and the noise. 
It reads, in Stratonovich form, 
\begin{align}\label{sto-Ham-intro}
    \begin{split}
      \dd q_i &= \frac{\partial H}{\partial p_i} \dd t + \sum_{l=1}^J\sigma_l(q_i) \circ \dd W_t^l
    \,,\\
    \dd p_i &= -\frac{\partial H}{\partial q_i} \dd t 
    - \sum_{l=1}^J \frac{\partial}{\partial q_i}\left (p_i \cdot\sigma_l(q_i)\right )  \circ  \dd W_t^l \, .
    \end{split}
\end{align}
We will refer to this system as the {\it Eulerian model}.

\section{Bridge dynamics}\label{sec:bridge}

All  models of the previous section can be written in the general It\^{o} form 
\begin{equation}\label{eq:x} 
\dd X_t = b(t,X_t) \dd t + \si(t,X_t) \dd W_t,\qquad  X_0 =x_0, \qquad t\in [0,T]\, , 
\end{equation}
where $b\colon [0,T]\times  \RR^N \to \RR^N$ and  $\si\colon [0,T] \times \RR^N \to \RR^{N\times N'}$ are  the drift and diffusion coefficient respectively. $W_t$ is a Wiener process in dimension $\RR^{N'}$ with independent components, and the stochastic integrals are It\^o.

In its most basic form, the previous problems correspond to simulating trajectories of $X:=(X_t,\, t\in [0,T])$, where both $X_0$ and $X_T$ are partially observed. 
We will encode the partial observation with two possibly different  projection operators $L_0$ and $L_T$, such that only the vectors $L_0X_0$ and $L_TX_0$ are observed.
As an example, if only landmarks positions are observed at times $0$ and $T$, then $L_0=L_T$ and $L_0$ is such that $L_0 X_t$ extracts only those components of the vector $X_t$ that correspond to landmark positions, thus ignoring the momentum components.

To model {\it inexact matching}, the observed vectors have extrinsic noise, and the realisations of random vectors
\begin{equation}\label{eq:conditions_V0_VT} 
    V_0 \sim N(L_0 X_0, \Sigma_0) \qquad \text{and} \qquad V_T \sim N(L_T X_T, \Sigma_T)
\end{equation}
are observed in instead of $L_0X_0$ and $L_TX_T$.
Notice that the {\it exact matching} limit is obtained upon taking $\Sigma = \eps^2 I$ and $\eps \rightarrow 0$. 
 
Simulating $X$ conditional on $(V_0, V_T)$ is challenging for several reasons:
\begin{enumerate}
  \item Contrary to many papers on diffusion bridge simulations (e.g.\ \cite{DelyonHu,beskos_exact_2006,PapaspiliopoulosRoberts,bladt_simulation_2016,schauer_guided_2017}), the conditioning is not on the full state, but on a subspace given by $L_0$ and $L_T$;
  \item the diffusion coefficient $\si$ can be state-dependent, such as in the Eulerian model;
  \item the dimension of the driving Brownian motion $W$ may not equal the dimension of the state space of the diffusion; and
  \item the regularity of the paths of landmark locations and momenta may not be the same (hypo-ellipticity), such as in the Lagrangian model.
\end{enumerate}

We demonstrate how these challenges can be tackled with the approach of~\cite{schauer_guided_2017,vdm-s-estpaper,partialbridges,smoothing_new}, where new methods for simulating conditioned diffusions using {\it guided proposals} were presented.
In addition to these works, the current setting needs special attention to obtain efficient numerical methods, due to the high dimension of the state space for large number of landmarks (given as $N=2dn$, where $d$ and $n$ denote the number and dimension of landmarks landmarks respectively), and to the nonlinear nature of the problem (mostly in choice of parameters and a particular "auxiliary process", described in Section \ref{sec:guided_ms_ahs}). Also, the structure of the landmarks configuration makes naive simple random-walk type Metropolis-Hastings schemes practically useless. As we will show, adaptations that take the geometry  properly into account, such as the  Riemannian Adjusted Metropolis Adjusted Langevin Algorithm, can be exploited instead.

\subsection{Conditioned diffusions} 

Throughout, we will assume that $X$ admits smooth transition densities $p$, such that $\PP(X_t \in \dd y \mid X_s =x) = p(s,x;t,\dd y)$ for $s<t$. In this section, we consider the problem of simulating the diffusion process $X$  conditioned on
\begin{equation}\label{eq:conditionings_x0vT} 
    X_0=x_0 \qquad \text{and}\qquad  V_T=v_T\, ,
\end{equation}
where $V_T \sim N(L_TX_T, \Sigma)$ and $x_0, v_T$ are given vectors.
This condition is different from \eqref{eq:conditions_V0_VT} as here, the initial state is fully observed, which turns out to be an important step in the algorithms that we present. 

The  diffusion conditioned on the events specified in \eqref{eq:conditionings_x0vT} is rather loosely referred to as a bridge process, as it bridges available observations on the state of the process at times $0$ and $T$. 

It is known that bridge processes satisfy an SDE similar to that of the unconditioned process, but with an additional {\it guiding} term in the drift.  
This additional drift term can be obtained  using Doob's $h$-transform (\cite{RogersWilliams2} or the theory of enlargement of initial filtration (\cite{JeulinYor,Baudoin})) to yield the bridge process  $X^\star=(X^\star_t,\, t\in [0,T])$ as the strong solution of the SDE
\begin{equation}\label{eq:xstar} \dd X^\star_t =b(t,X^\star_t) \dd t + a(t,X^\star_t) r(t,X^\star_t) \dd t + \si(t,X^\star_t) \dd W_t, \qquad X^\star_0=x_0, \end{equation}
where $a=\si \si^\T$ (where $^\T$ denotes matrix transpose), $r(t,x)=\nabla_x \log \rho(t,x)$, $t\in [0,T)$, $x\in \RR^N$ and
\begin{equation}\label{eq:inexactrho}
	\rho(t,x) = \int_{\RR^N} p(t,x; T,\xi) \psi(v_T; L_T\xi, \Sigma) \dd \xi\, .
\end{equation}
Here, $\psi(x;\mu,\Sigma)$ denotes the density of the $N(\mu,\Sigma)$-distribution, evaluated at $x$, for any $\mu\in \mathbb R^N$ and $\Sigma \in \mathbb R^{N\times N}$. Whereas, we focus here on inexact matching, we remark that 
in case of exact matching (where $\Sigma \equiv 0$) the formula for $\rho$ is more complex, see for example Section 1.3.2 in \cite{partialbridges}.
\begin{ex}
To help gain the reader some intuition on the term $a(t,X^\star_t) r(t,X^\star_t)$, consider the setting where $b\equiv 0$, $\sigma$ is constant and the diffusion is fully observed at time $T$ with $N(0,\eps^2 I)$-noise. In that case, the SDE for the conditioned process is given by $\dd X^\star_t = (T+\eps^2-t)^{-1} (v_T-X^\star_t) \dd t + \sigma \dd W_t$ (this follows from first computing $\rho$ in \eqref{eq:inexactrho} and subsequently $r$).  Note that for small values of $\eps$, the guiding term is roughly ``remaining distance to endpoint'' divided by ``remaining time'' and that a value of $\epsilon$ that is strictly positive avoids dividing by zero in the drift at time $T$. 
\end{ex}

\subsection{Bayesian formulation  for landmarks matching}
We will follow the Bayesian paradigm for statistical inference in this paper. 
Within Bayesian statistics it is common to write the data-generating model in a hierarchical way. This hierarchical scheme (to be read from bottom upwards) is as follows:
\begin{gather*}
	v_T \mid x_T  \sim N(L_T x_T,  \Sigma_T) \\
	v_0 \mid x_0  \sim  N(L_0 x_0, \Sigma_0) \\
	x_T \mid x_0  \sim p(0,x_0; T, x_T) \\
	x_0  \sim p(x_0)\, .
\end{gather*}
Here, $p(x_0)$ denotes the prior density on the initial state.  The observed variables are $v_0$ and $v_T$; the unobserved variables are $x_0$ and $x_T$. Bayesian inference is based on the distribution of unobserved variables, conditional on observed variables. Hence, using Bayesian notation
\[
	p(x_0, x_T \mid v_0, v_T)  \propto p(v_T \mid x_T) p(v_0 \mid x_0) p(x_T\mid x_0) p(x_0) 
\]
(note that the third term on the right-hand-side is the intractable transition density). 
In particular, this implies that 
\begin{equation}\label{eq:fullcond_xT} p(x_T \mid v_0, v_T, x_0) \propto p(v_T \mid x_T) p(x_T \mid x_0). \end{equation}
As shown in  the appendix of~\cite{smoothing_new}, the solution of equation \eqref{eq:xstar} with \eqref{eq:inexactrho} implies that $X^\star_T$ has density 
\begin{equation}\label{eq:eta}
	\xi \mapsto \eta(\xi; x_0, v_T) := \frac{p(0,x_0; T,\xi) \psi(v_T;L_T\xi, \Sigma) \dd \xi}{\int p(0,x_0; T, \xi')  \psi(v_T; L_T\xi', \Sigma) \dd \xi'}, 
\end{equation}
but this is just the density in \eqref{eq:fullcond_xT}. This reveals that  sampling of $x_T$, conditional on $(v_0, x_0, v_T)$, can be done by sampling  $X^\star_T$, which in turn is obtained as the endpoint of a sampled bridge path.

\section{Guided proposals}
\subsection{Overview of  guided proposals}
Since the transition densities of a diffusion are only available in closed form in very special cases, $\rho$ in \eqref{eq:inexactrho} is intractable. This prevents using a numerical discretisation scheme on the SDE for $X^\star$, where the process would be simulated forward on a fine grid (for example using Euler's method). The key idea of~\cite{schauer_guided_2017} consists of instead simulating from a process $X^\circ$ that is obtained after replacing $p$ by the transition densities $\widetilde{p}$ of an auxiliary diffusion process $\widetilde{X}$ for which $\widetilde{p}$ is tractable.  
This is naturally the case for linear processes $\widetilde{X}$ whose dynamics are governed by the SDE
\begin{equation}\label{eq:auxlin} 
    \dd \widetilde{X}_t =\widetilde{b}(t,\widetilde{X}_t) \dd t + \widetilde\si(t) \dd W_t\, , 
\end{equation}
where $\widetilde b$ is of the form
\begin{align}
	\widetilde{b}(t,x) = \widetilde\beta(t) +\widetilde{B}(t) x\, . 
\end{align}
Hence, instead of sampling from $X^\star$, one samples from $X^\circ$ defined by
\begin{equation}\label{eq:xcirc} 
    \dd X^\circ_t =b(t,X^\circ_t) \dd t + a(t,X^\circ_t) \widetilde{r}(t,X^\circ_t) \dd t + \si(t,X^\circ_t) \dd W_t, \qquad X^\circ_0=x_0\, ,
    \end{equation}
with $\widetilde{r}(t,x)=\nabla_x \log \widetilde\rho(t,x)$, where $\widetilde\rho$ is derived from $\widetilde{p}$ exactly as $\rho$ is in terms of $p$. We call  $t\mapsto a(t,X^\circ_t) \widetilde{r}(t,X^\circ_t)$ the {\it guiding term}, as it is a term that is superimposed on the drift of the original SDE for $X$ to satisfy the imposed conditioning at time $T$.
Draws from this {\it guided diffusion bridge proposal} can subsequently be accepted/rejected in a Metropolis-Hastings sampler, provided that 
\begin{enumerate}
  \item the law of $X^\star$ is absolutely continuous with respect to the law of $X^\circ$;
  \item the Radon-Nikodym derivative showing up in the acceptance probability is tractable.
\end{enumerate}
A short introductory account on the Metropolis-Hastings algorithm is given in Section \ref{sec:mh}. 
\begin{table}[h]
\begin{tabular}{|l|l|l|l|l|}
\hline
$X$ & original, unconditioned diffusion process, defined by \eqref{eq:x} & $b$ & $\si$& $\PP$\\
$X^\star$ & corresponding  bridge, conditioned on $v$, defined by \eqref{eq:xstar}& $b^\star$ &$\si$& $\PP^\star$\\
$X^\circ$ & proposal process defined by \eqref{eq:xcirc}& $b^\circ$&$\si$& $\PP^\circ$\\ 
$\widetilde X$ & linear process defined by  \eqref{eq:auxlin} with transition densities $\widetilde p$ & $\widetilde b$ &$\widetilde\si$& $\widetilde{\PP}$\\ &
appearing in the definition of $X^\circ$& &&\\
\hline
\end{tabular}
\caption{Summary of the notation of various bridges in this work. The rightmost three columns give the drift, diffusion coefficient and measure on $C([0,T], \RR^N)$ respectively. }
\end{table}

An implicit assumption made throughout is that a strong solution to the SDE \eqref{eq:xcirc} exists, but, in addition, guided proposals only exist if the auxiliary process is chosen such that $\widetilde{X}$ admits transition densities, see Section 2.1 in \cite{partialbridges} for more details. 
Furthermore, a key result from~\cite{partialbridges} provides sufficient conditions for absolute continuity of $\PP^\star$ with respect to $\PP^\circ$, denoted as $\PP^\star \ll \PP^\circ$, and an expression for the Radon-Nikodym derivative
\begin{equation}\label{eq:LR}
    \frac{\dd \PP^\star}{\dd \PP^\circ}(X^\circ) = \frac{\widetilde \rho(0+,x_0)}{ \rho(0+,x_0)}\Psi(X^\circ)\, , 
\end{equation}
with $\rho$ defined  in \eqref{eq:inexactrho} and
\begin{equation}\label{eq:Psi} 
    \Psi(X^\circ)=\exp\left(\int_0^T \scr{G}(s,X^\circ_s) \dd s\right)\, , 
\end{equation}
where
\begin{equation}
\begin{split}\label{eq:G} \scr{G}(s,x) &= (b(s,x) - \widetilde b(s,x))^\T \widetilde r(s,x)  \\ & \qquad -  \frac12 \mbox{tr}\left(\left[a(s,x) - \widetilde a(s)\right] \left[\widetilde H(s)-\widetilde{r}(s,x)\widetilde{r}(s,x)^\T\right]\right)\, 
\end{split}
\end{equation}
and $\widetilde{H}(s)$ is the negative of the Hessian matrix of $x\mapsto \nabla_x \log \widetilde\rho(s,x)$, which turns out to be a constant on $x$. 

\medskip

Forward simulation of the guided proposal $X^\circ$ in \eqref{eq:xcirc} requires evaluation of $\tilde r$. It turns out that in all algorithms of Section \ref{sec:mcmc} the intractable term $\rho(0+,x_0)$, appearing in \eqref{eq:LR} and depending on the transition densities of the process $X$, cancels.  Ignoring this term, evaluating \eqref{eq:LR}  requires evaluation of $\tilde r$, $\tilde H$ and $\tilde\rho(0+,x_0)$. These quantities can be computed efficiently as they only depend on the tractable process $\tilde X$. There exist various expressions for this, full details are in \cite{smoothing_new}. For completeness we recap  Theorem 2.4 in \cite{smoothing_new} which gives for $t\in [0,T]$
\begin{align*} \tilde r(t,x) &= L(t)^\T M(t)  (v(t) -\mu(t) - L(t) x) \\
\tilde{\rho}(t,x) &= \psi(v(t); \mu(t)+L(t) x, M^\dagger(t))	\\
H(t) &=L(t)^\T M(t) L(t)
\end{align*}
Here $M(t) = [M^\dagger(t)]^{-1}$ with $(L(t), \mu(t), M^\dagger(t))$ being defined as solutions to the (backward) differential equations 
 \begin{equation}\label{eq:ode}
 \begin{split}	
 \dd L(t) = -L(t) \tilde B(t)\dd t \qquad & L(T)=L_T \\
 \dd \Md(t)=- L(t) \tilde{a}(t) L(t)^\T\dd t \qquad & M^\dagger(T)=\Sigma_T  \\
 \dd \mu(t) =-L(t) \tilde\beta(t)\dd t \qquad & \mu(T)=0 
 	 \end{split}
 \end{equation}
and $v(t)=v_T$ for $t\in (0,T]$. 
 
As we assume  that $V_0 \sim N(L_0 X_0, \Sigma_0)$ (as in \eqref{eq:conditions_V0_VT}), it follows again from Theorem  2.4 in \cite{smoothing_new} that 
 for $t=0$
  \begin{equation}\label{eq:change-at-S} L(0)=\Bm L_0 \\ L(0+) \Em, \quad 
 \Md(0)=\Bm \Sigma_0 & 0 \\ 0 & M^\dagger(0+)\Em,  \quad 
 \mu(0)=\Bm 0 \\ \mu(0+) \Em, \quad v(0)= \Bm v_0 \\ v_T \Em 
\end{equation}
(with $0$ denoting a vector or matrix of appropriate dimensions). 
Note that for fixed $x$, $t\mapsto \tilde\rho(t,x)$ is left-continuous.

\begin{rem}
From the discussion in Appendix \ref{app:matching_conditions} on the exact matching case, it is in general preferable to take $\widetilde b$ and $\widetilde \sigma$ such that 
\begin{equation}\label{eq:match-conject} 
    L_T b(T,X^\circ_T) = L_T \widetilde{b}(T,X^\circ_T) \quad \text{and} \quad  L_T a(T,X^\circ_T) L_T^\T = L_Ta(T)L_T^\T\, ,
\end{equation}
which is possible for both the Lagrangian  and Eulerian model, as we show in the next section. 
\end{rem}

\subsection{Guided proposals for  landmark models}\label{sec:guided_ms_ahs}

Guided proposals are specified by choice of the auxiliary process $\widetilde{X}$. Here we detail this choice for the three landmarks models considered. Recall (see also Equation \eqref{hamiltonian}) 
\begin{equation}\label{eq:deriv-deterministic} 
    \frac{\partial H}{\partial p_i} = \sum_{j=1}^n p_j k(q_i-q_j) \qquad \text{and} \qquad \frac{\partial H}{\partial q_i}=\sum_{j=1}^n \langle p_i, p_j\rangle \nabla k(q_i-q_j)\, . 
\end{equation}

\subsubsection{Choice of auxiliary process for the Langevin  and Lagrangian model}

In the Langevin model we have
\begin{align*}
    \dd  q_i^\alpha &= \frac{\partial H}{\partial p_i^\alpha} \dd t \qquad \mathrm{and} \qquad 
	\dd p_i^\alpha =-\lambda \frac{\partial H}{\partial p_i^\alpha} \dd t   -\frac{\partial H}{\partial q_i^\alpha} \dd t + \gamma_i \dd W_t^i\, .
\end{align*}
The matching conditions \eqref{eq:match-conject} suggest to match both the drift on the landmarks positions and the diffusivity at time $T$. 
In view of \eqref{eq:deriv-deterministic}, we take the  auxiliary process $\widetilde{X}$ as 
\begin{align*}
    \dd  \widetilde{q}_i^\alpha &= \sum_{j=1}^n \widetilde{p}_j k(q^T_i-q^T_j) \dd t \qquad \mathrm{and} \qquad 
	\dd p_i^\alpha =-\lambda \sum_{j=1}^n \widetilde{p}_j k(q^T_i-q^T_j)\dd t   + \gamma_i \dd W_t^i\, ,
\end{align*}
where $q^T_i$ is the observed $i$-th landmark position at time $T$. Because the kernels are evaluated on the final positions, these equations define a linear process, from which $\widetilde{B}$, $\widetilde\beta$ and $\widetilde{\si}$ can be inferred directly. 

\subsubsection{Choice of auxiliary process for the Eulerian model}\label{subsec:choice_aux}

As the matching assumptions are formulated for an SDE in It\^{o} form, we first need to transform the Eulerian model equations from Stratonovich to It\^{o} form. To compute the additional term in the drift, we specify the noise kernels to be located at fixed positions. We suppose  at location $\delta_\ell$ kernel $q\mapsto \si_\ell(q)$ with coordinates as specified in Equation \eqref{eq:noisefield}.
\begin{prop}\label{lem:ahs-strat2ito}
For the Eulerian model,	the additional term in the drift when switching from Stratonovich to It\^{o} form for $q$ is given by 
\begin{equation}
	\label{eq:stratito-q}
	 \frac12 \sum_\ell z_\ell(q) \bar{k}_\tau(q-\delta_\ell) \ga\, .
\end{equation}
For $p$ this term is given by
\begin{equation}
	\label{eq:stratito-p}
 \frac12 \sum_\ell \langle p,\ga\rangle \left( z_\ell(q) \nabla \bar{k}_\tau(q-\delta_\ell)   - \bar{k}_\tau(q-\delta_\ell) \nabla z_\ell(q)\right)\, .  
\end{equation}
Here, in the notation we have  omitted dependence of $(\tau, \delta, \ga)$ on $\ell$ and have denoted 
\begin{align*}
    z_\ell(q) = \langle \nabla \bar{k}_\tau(q-\delta_\ell), \ga\rangle\,.
\end{align*}

\end{prop}

This proposition shows that there is a simple way to take the Stratonovich-It\^{o} correction into account in the drift of the auxiliary process $\widetilde{X}$: for landmark $i$ one superimposes the terms \eqref{eq:stratito-q} and \eqref{eq:stratito-p} to the drift of the auxiliary process for the Lagrangian model, with the final positions $q=q_i^T$ substituted.

\section{MCMC algorithms for landmark matching and template estimation}\label{sec:mcmc}

Both the drift $b$ and the diffusion coefficient $\si$ in \eqref{eq:x} may contain unknown parameters. Suppose the vector of unknown parameters is $\th$ and we will write for example $b_\th$ and $\si_\th$ to highlight the explicit dependence on these parameters.  
 In this section, we present an MCMC (Markov Chain Monte Carlo) algorithm for landmarks matching (Section \ref{subsec:lm_matching}) and template estimation (Section \ref{sec:rmmala-steps}), including estimation of unknown parameters. 

\subsection{Bridging two landmark configurations}\label{subsec:lm_matching}

The Bayesian approach to the problem of landmarks matching boils down to sampling from $(X, \th)$ conditional on $(v_0, v_T)$ (the observed two landmarks configurations). It is well known that in case of unknown parameters in the diffusion coefficient, a Gibbs sampler that successively updates $\th \mid (X, v_0, v_T)$ and $X \mid (\th, v_0, v_T)$ lead to an invalid (reducible) scheme (see for example~\cite{RobertsStramer}). To circumvent this problem, we will update the Wiener increments $W$ instead of the process $X$ itself.
 Then, as we have made the assumption that  a strong solution to the SDE \eqref{eq:xcirc} exists, there is a measurable map $\GP_\th$ such that $X^\circ=\GP_\th(x_0, W)$, where $W$ is the driving Wiener process in $\RR^{N'}$ ($\GP$ being an abbreviation of $\scr{G}$uided $\scr{P}$roposal). 	Note that we could have written $\GP_\th(x_0, W, v_T)$ instead of $\GP_\th(x_0, W)$ to highlight that the guided proposal always depends on the state observed with error at time $T$. To reduce notational overhead we have chosen not to do so. 

Let $x_0=(q_0, p_0)$ and $x_T=(q_T, p_T)$ denote the latent states.
We propose a Gibbs sampler for updating $(W, \theta, p_0)$ and assume that the noise on $v_0$ is negligible so that it can be assumed that $q_0=v_0$. Then, for each sample of $(W, \theta, p_0)$ we obtain a sample of the guided proposal via $\GP_\theta((q_0, p_0), W)$, which henceforth also yields a sample of $x_T$. 

After initialisation of  $(W, \theta, p_0)$  the sampler cycles over the  steps
\begin{enumerate}
  \item sample  $W$ conditional $(x_0, \th, v_T)$ using guided proposals, with pCN (preconditioned Crank-Nicolson) updates on the Wiener increments with Algorithm \ref{alg1};
  \item sample  $p_0$ conditional on $(q_0, \th,  W, v_T)$ using MALA (Metropolis Adjusted Langevin Algorithm), with gradients obtained using automatic differentiation with Algorithm \ref{alg2};
  \item sample $\th$ conditional on $(x_0, W, v_T)$ with Algorithm \ref{alg3}. 
\end{enumerate}
In these algorithms, we denote  prior densities  on $\th$, $q_0$ and  $p_0$ generically by $\pi$ (Bayesian notation). In each of the steps, we keep track of the changes to $X$ induced by either updating $W$, $p_0$ or $\theta$. The acceptance probabilities in Algorithms \ref{alg1}, \ref{alg2}, \ref{alg3} and \ref{alg:rmmala} follow from the computations in Section 4.1 of \cite{smoothing_new}.   For readers less familiar with Markov Chain Monte Carlo algorithms, we refer to Section \ref{sec:mh}, where we also derive the acceptance probability for a simplified version of Algorithm \ref{alg1}.

\begin{algorithm}[h]
 \begin{algorithmic}[1]
\STATE{ Choose a persistence parameter $\eta \in [0,1]$}
  \STATE {Compute $L_\th(t)$, $\Md_\th(t)$ and $\mu_\th(t)$ for $t\in [0,T]$ by solving \eqref{eq:ode}.    \item Sample a Wiener  process $Z$, independently of $W$ and set 
  \begin{align*}
    W^\circ=\eta W + \sqrt{1-\eta^2} Z\, .
  \end{align*}}
 \STATE  Compute  $X^\circ=\GP_\th(x_0, W^\circ)$.
\STATE Compute 
\begin{align*}  
    A = \Psi_\th(X^\circ)/\Psi_\th(X)\, .
\end{align*}
\STATE Draw $U\sim \scr{U}(0,1)$.
  \STATE if $U<A$ 
set  $X=X^\circ$ and $W=W^\circ$\;
 else
   keep previous $X$ and $W$.
 \caption{Update $(W, X)$, conditional on $(x_0, \theta, v_T)$}
   \label{alg1}
  \end{algorithmic}
\end{algorithm}

 \begin{algorithm}[h]
 \begin{algorithmic}[1]
\STATE Choose a step-size $\delta>0$\;
\STATE  Sample $Z \sim N_{dn}(0, I)$ (with $d$ the dimension of a landmark). 
\STATE { Set $x_0 =(q_0, p_0)$ and  propose
\begin{align*}  
p_0^\circ = p_0 + \frac{\delta}{2}  \scr{L}_\th(x_0,W) + \sqrt{\delta} Z\, ,
\end{align*}
where 
\begin{align*}  
\scr{L}_\th(x_0,W) = \nabla_{p_0}\Big( \log \Psi_\th(\GP_\th(x_0, W)) + \log \widetilde\rho_\th(0,x_0)\Big)\, . 
\end{align*}
Here, the gradient is obtained by automatic differentiation. }
\STATE { Set $x_0^\circ   =(q_0,  p_0^\circ)$ and compute $X^\circ = \GP_\th(x_0^\circ,W)$ and $\scr{L}_\th(x_0^\circ,W)$. 
	\item Compute 
\begin{align*}  
    A = \frac{\Psi_\th(X^\circ)}{\Psi_\th(X)} \frac{\widetilde\rho_\th(0,x^\circ_0)}{\widetilde\rho_\th(0,x_0)}  \frac{\pi(p_0^\circ)}{\pi(p_0)}\frac{\psi(p_0; p_0^\circ + \delta \scr{L}_\th(x_0^\circ,W)/2 , \delta I)}{\psi(p_0^\circ; p_0 + \delta \scr{L}_\th(x_0,W)/2 , \delta I)}\, .
\end{align*}}
\STATE Draw $U\sim \scr{U}(0,1)$.

\STATE if $U<A$, set  $X=X^\circ$ and $p_0=p_0^\circ$ else  keep previous $X$ and $p_0$

 \caption{Update $(p_0, X)$, conditional on $(q_0,\th,  W, v_T)$}
 \label{alg2}
\end{algorithmic}
\end{algorithm}

 \begin{algorithm}[h]
   \begin{algorithmic}[1]
\STATE Choose an irreducible Markov kernel $q$ on the domain of $\th$
\STATE Sample $\th^\circ$ from a kernel $q(\cdot \mid \th)$. 
\STATE Compute $X^\circ = \GP_{\th^\circ}(x_0, W)$. 
\STATE{ Compute 
\begin{align*}
    A = \frac{\Psi_{\theta^\circ}(X^\circ)}{\Psi_\th(X)} \frac{\widetilde\rho_{\theta^\circ}(0,x_0)}{\widetilde\rho_\theta(0,x_0)} \frac{\pi(\th^\circ)}{\pi(\th)} \frac{q(\th \mid \th^\circ)}{q(\th^\circ \mid \th)}\, .
\end{align*}}
\STATE Draw $U\sim \scr{U}(0,1)$.
\STATE if $U<A$,  set $X=X^\circ$ and  $\th=\th^\circ$ else  keep previous $X$ and $\th$

 \caption{Update $(\th, X)$, conditional on $(x_0, W, v_T)$}
 \label{alg3}
 \end{algorithmic}
\end{algorithm}

\subsection{Template estimation}\label{sec:rmmala-steps}

In the template estimation problem, we assume $I$ landmark configuration are observed with error. We assume that each configuration is a stochastic deformation of a fixed template configuration $q_0$. Fix an initial momentum vector $p_0$ and denote  $x_0=(q_0,p_0)$.  Let  $\{W^i,\, i=1,\ldots, I\}$ be independent Wiener processes on $[0,T]$.  We assume that the $i$-th  landmark configurations is a sample from $v^i_T = L_T X^i_T + N(0, \Sigma_T)$, where  $X^i$ has dynamics governed by the SDE of the Lagrangian , Langevin , or Eulerian model with initial state $x_0$ and driven by the Wiener process $W^i$. 

For simplicity we will assume the initial momentum to be zero, though it would not pose a problem to relax this assumption using Algorithm \ref{alg2}. A Gibbs sampler to sample from the posterior distribution of $(W^1,\ldots, W^I, \theta, q_0)$ is obtained by iterating the steps
\begin{enumerate}
	\item sample each $W^i$ conditional on $(x_0,\theta, v^i_T)$ independently using Algorithm \ref{alg1};
	\item sample $\theta$ conditional on $(x_0, W^1,\ldots, W^I, v^1_T,\ldots, v^I_T)$ using a slight adaptation of Algorithm \ref{alg3}, where the loglikelihood contributions from different $W^i$-s are added;
	\item sample $q_0$ conditional on $(W^1,\ldots, W^I, v^1_T,\ldots, v^I_T)$.
\end{enumerate}
For updating $q_0$, with a large number of landmarks simple random-walk updates will perform terribly bad    as these do not respect the geometry of the landmark configuration. For that reason we propose to  use the RMMALA (Riemannian Manifold Metropolis Adjusted Langevin Algorithm)   as introduced in Section 5 of~\cite{girolami-calderhead}. Details are given in Algorithm \ref{alg:rmmala}.

  \begin{algorithm}[h]
  \begin{algorithmic}[1]
\STATE  Choose a step-size $\delta>0$\;
\STATE Sample $Z \sim N_{dn}(0, K(q_0))$, where  $K(q)$ is the matrix with blocks $K(q)_{ij}=K(q_i, q_j)$, $1\le i \le n$.
\STATE Propose 
\begin{align*}
 q_0^\circ = q_0 + \frac{\bar\delta}{2} K(q_0) \scr{L}(q_0) + \sqrt{\delta}Z\, ,
\end{align*}
 where 
\begin{align*}
\scr{L}(q_0) :=  \nabla_{q}\Big( \log \widetilde\rho_\th(0,(q_0,p_0)) + \sum_{i=1}^I \log \Psi_\th(g((q_0,p_0), W^i)) \Big)\, .
\end{align*}
Here, the gradient is obtained by automatic differentiation. Set $x_0^\circ =(q_0^\circ, p_0)$. 
\STATE For $i=1,\ldots, I$, compute $X^{\circ,i} = \GP_\th(x_0^\circ,W^i)$ and $\scr{L}_\th(q_0^\circ)$. 
\STATE Compute 
\begin{align*}
    A =  \frac{\widetilde\rho_\th(0,x_0^\circ)}{\widetilde\rho_\th(0,x_0)}  \frac{\pi(q_0^\circ)}{\pi(q_0)} \frac{\psi(q_0; q_0^\circ + \delta K(q_0^\circ)\scr{L}(q_0^\circ)/2 , \delta K(q_0^\circ))}{\psi(q_0^\circ; q_0 + \delta K(q_0) \scr{L}(q_0)/2 , \delta K(q_0))}  \prod_{i=1}^I \frac{\Psi_\th(X^{\circ,i})}{\Psi_\th(X^i)}\, . 
\end{align*}
\STATE Draw $U\sim \scr{U}(0,1)$.
\STATE if $U<A$, set  $X^i=X^{\circ, i}$ for $i=1,\ldots, I$ and $q_0=q_0^\circ$ else keep previous $X$ and $q_0$
  \caption{Update $q_0$, conditional on $(\theta, W^1,\ldots, W^I, v_T^1,\ldots, v_T^I)$}
\label{alg:rmmala}
    \end{algorithmic}

\end{algorithm}

\subsection{Computational aspects}
The computational cost of the proposed algorithms manifests itself in two ways. First, if parameter estimation is included (Algorithm \ref{alg3}), then the backward ODEs in Equation \eqref{eq:ode} have to be recomputed  if the parameters of the proposal change. Note however that  in state-of-the-art competing methods the parameter is fixed and if we choose to do likewise, then  the ODEs need to computed only once. Once solved on a fine temporal grid,
at each gridpoint the Cholesky decomposition of $\Md_t$ needs to be computed, which scales cubically in the number of landmarks. 
With a large number of landmarks the computational costs can drastically be decreased by either  applying an ensemble backward filter or enforcing sparsity. For details of these approaches we refer to Sections 6 and 7.3 in \cite{smoothing_new}.  

Secondly, both Algorithms \ref{alg2} and \ref{alg:rmmala} rely on gradients obtained by 
automatic differentiation. The runtime of the algorithms therefore 
depends on the computational efficiency of the automatic differentiation 
algorithm used. With large numbers of landmarks we expect reverse-mode automatic differentiation to be more efficient. Note that this is an implementational issue and not a deficit of the proposed methods. Finally, we remark that the Eulerian model is computationally more demanding than the Langevin model, as its diffusivity is non-constant. This incurs some extra cost in simulating the guided proposal but, more importantly, also increased computational cost in computing $\Psi(X^\circ)$ (cf.\ Equation \eqref{eq:Psi}). Contrary to the Lagrangian model, for the Eulerian model also the second term on the right-hand-side of Equation \eqref{eq:G} needs to be evaluated.

From a practical point of view, we have noticed that tuning of stepsizes in the algorithm can be delicate. This can potentially be resolved by adaptive tuning at early iterations and fixing the stepsizes at subsequent iterations.

\section{Numerical examples}\label{sec:numexamples}

\subsection{Settings} Before presenting examples we discuss precise settings in the numerical experiments.

\subsubsection{Specification of Hamiltonian kernel}

We take a Gaussian kernel $K$ of the form $K(x) = k(x) \mbox{Id}_d$ with 
\begin{align}
    k_a(x) = c \exp\left( - \|x\|^2/(2a^2)\right)\, .
    \label{kernel}
\end{align}
Note that $\nabla k_a(x) = -c\, a^{-2}\, k(x)\, x$. Without loss of generality we  take $c=1$: 
\begin{rem}
For the Lagrangian model specified in Equation \eqref{TV-model}, we have, by writing $\tilde H = c H$, 
\begin{align}
    \dd  q_i^\alpha &= c \frac{\partial \tilde H}{\partial p_i^\alpha} \dd t \qquad \mathrm{and} \qquad \dd p_i^\alpha = -c \frac{\partial \tilde H}{\partial q_i^\alpha} \dd t + \gamma_i \dd W_t^i\, , 
\end{align}
which, upon rescaling $\gamma_i \rightarrow \gamma_i \sqrt{c}$ and $t \rightarrow t / c$ (and $T \rightarrow T / c$), has the form of \eqref{TV-model} with $c=1$. This implies that the parameters $\gamma$ and $c$ are not mutually independent. A similar consideration applies to the Eulerian model.
\end{rem}

\subsubsection{Specification of intrinsic  noise}\label{subsubsec:noisespec}
For the Lagrangian model, we take the noise on all landmarks the same, scaled by $\sqrt{n}$. That is,  $\gamma_i =\gamma/\sqrt{n}$ for all $i\in \{1,\ldots,n\}$. In the experiments we take $\gamma=0.1$. 

For the Eulerian model, let $\tau>0$. We take the noise fields centred at points that are both horizontally and vertically separated by a distance that is an integer multiple of $2\tau$. Denote the locations of the noise fields by $\{\delta_j\}$. At each location $\delta_j$ we take noise fields to be
\begin{align*}
\frac{2}{\pi}\gamma  k_\tau(x-\delta_j)\Bm 1 \\ 0 \Em \qquad \text{and} \qquad \frac{2}{\pi}\gamma   k_\tau(x-\delta_j)\Bm 0 \\ 1 \Em\, .
\end{align*}
 This kernel is the same as the Hamiltonian kernel, but does not have to be so in general. 
The scaling by $2/\pi$ is chosen so that if $\gamma=1$ the noise on the landmarks positions is close to $1$ uniformly in space (close to a decomposition of unity). 

\subsubsection{Specification of extrinsic noise}\label{subsec:eps}
In each example we choose $\Sigma_T=\epsilon^2 I$ with  $\epsilon=0.01$ (cf.\ Equation \eqref{eq:inexactrho}), except in the examples of sections \ref{subsec:1Dlandmarks} and \ref{subsec:cardiac} where we took $\eps=0.001$.
From a practical point of view, as landmarks are often annotated manually, it is a realistic assumption that $\epsilon$ is strictly positive.

\subsubsection{Specification of parameter $\th$}

In the simulations, we fix the scaling parameter $\tau$ of the noise-kernels, as well as their positions. We estimate the Hamiltonian kernel $a$, but fix  the size (amplitude) of the noise ($\gamma$). Hence, $\th=a$.

\subsubsection{Prior on the landmark positions and momenta}
In the problem of bridging two landmark configurations we use a prior on the landmark momenta similar to the prior suggested in Section 3 of~\cite{Marsland-Shardlow}. More specifically, we take for $\kappa_{\rm mom}, \kappa_{\rm pos}>0$, 
\begin{align*} 
    \pi(q_1,\ldots, q_n, p_1,\ldots, p_n) &= \pi(p_1,\ldots, p_n\mid q_1,\ldots, q_n) \pi(q_1,\ldots, q_n) \\ 
    &= \psi\left(p_1,\ldots, p_n; 0, \kappa_{\rm mom} K(q)^{-1}\right) \prod_{i=1}^n \psi(q_i; 0, \kappa_{\rm pos})\, .
\end{align*}
As we assumed that the landmark positions are observed, we make the shortcut where we assume $(p_1,\ldots, p_n) \sim N(0, \kappa_{\rm mom} K(q_0)^{-1})$,  Here, $\kappa_{\rm mom}$ is a parameter that we fix to a large value.  In the experiments we take $\kappa_{\rm mom}=100$. 

\subsubsection{Tuning parameters for the MCMC-sampler}

As can be seen from algorithms \ref{alg1}, \ref{alg2} and \ref{alg3}, the proposed MCMC-algorithm requires specification of 
\begin{itemize}
  \item the persistence parameter $\eta \in [0,1]$ (Algorithm \ref{alg1});
  \item the step size $\delta >0$ (Algorithm \ref{alg2});
  \item the Markov kernel $q$ for updating $\th$ (Algorithm \ref{alg3}).
\end{itemize}
Naturally, these choices affect the efficiency of the algorithm to explore the support of the target distribution. We target acceptance rates of about $50\%$. For the Markov kernel $q$ as $\th$ we use updates of the form $\log \th_i^\circ \mid \th_i \sim N(\log \th_i, \si^2_{\th})$ to ensure $\th^\circ_i>0$.

\subsubsection{Time discretisation}\label{sec:timediscretisation}
We discretise guided proposals in time by taking a regular grid on $[0,1]$ with mesh-width $0.01$, followed by applying the map $s\mapsto s(2-s)$ to this grid. The latter ensures the discretisation is finer near the endpoint $T$, where the guiding term can become large  (see also Section 5 on~\cite{vdm-s-estpaper} for discretisation of guided proposals, where the choice of maps $s\mapsto s(2-s)$ is motivated). 

For computing the guiding term and the right-hand-side of \eqref{eq:LR}, we solve the system in \eqref{eq:ode}, where the differential equation for $L(t)$ was solved with an implicit Euler scheme and both $\Md$ and $\mu$ were approximated using the trapezoid rule (all on a regular grid with mesh-width $0.01$).   

\subsubsection{Software -- hardware}
The source code of the numerical examples is available together with a Julia package \texttt{BridgeLandmarks} for shape analysis with stochastic landmark dynamics, \cite{BridgeLandmarks}. Automatic differentiation was implemented via \cite{RevelsLubinPapamarkou2016}.
All computations were done on a  MacBook Pro, with a 2 GHz Quad-Core Intel Core i5 with 16 GB RAM. 

\subsection{Examples for landmarks matching}

\subsubsection{1D landmarks}\label{subsec:1Dlandmarks} 
We start with an example that is meant to illustrate the developed methodology in a 1-dimensional setup. We consider a problem with just $3$ landmarks.
Here, the initial locations of the points are given by $-0.5$, $0.0$ and $0.1$ and we condition on $-0.5$, $0.2$ and $1.0$ respectively. We consider both the Lagrangian  and Eulerian model and fixed the Hamiltonian kernel parameter to $1.0$. For the Eulerian model, $6$ noise sources were positioned on the domain $[-2.5, 2.5]$ with $\tau=0.5$. For the Lagrangian model we fixed  $\gamma=1/\sqrt{n}\approx 0.58$, for the Eulerian model we took $\gamma=0.1$. Hence, we assume all parameters are known, except for the initial momenta and bridges. We consider landmarks matching without noise, which is forced by assuming that the final position is observed with additive $N(0, 10^{-6})$-noise ($\eps=0.001$, cf.\ Section \ref{subsec:eps}). 

We ran an algorithm composed of alternating one step of Algorithm \ref{alg1} (bridge updating) and one step of Algorithm \ref{alg2} (initial momenta updating) for $20\_000$ iterations. We computed the guided proposal on a grid with mesh-width $0.001$, subsequently mapped to a nonequidistant grid as explained in Subsection \ref{sec:timediscretisation}. We  initialised with zero momentum for all of the three points at time $0$.

\begin{figure}
\begin{center}
	\includegraphics[scale=0.6]{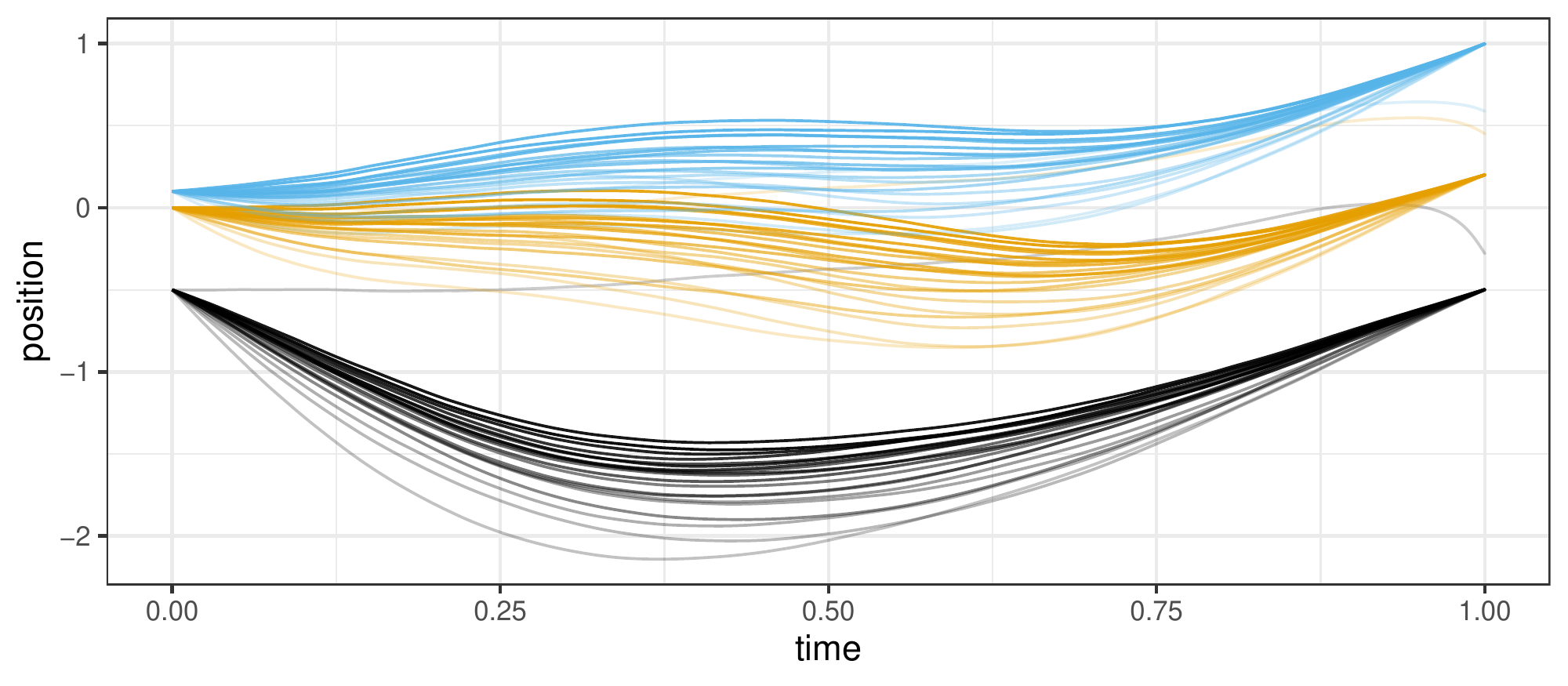}
\caption{Bridges for the 1D landmarks example, Lagrangian model.  The sampler ran for  $20\_000$ iterations. Every $1000$-th iterations is shown. ``Early'' iterations are depicted less thick. Note that the initial trajectory does not satisfy the conditioning.  \label{fig:bridges_1D_ms}}
\end{center}
\end{figure}

\begin{figure}
\begin{center}
	\includegraphics[scale=0.6]{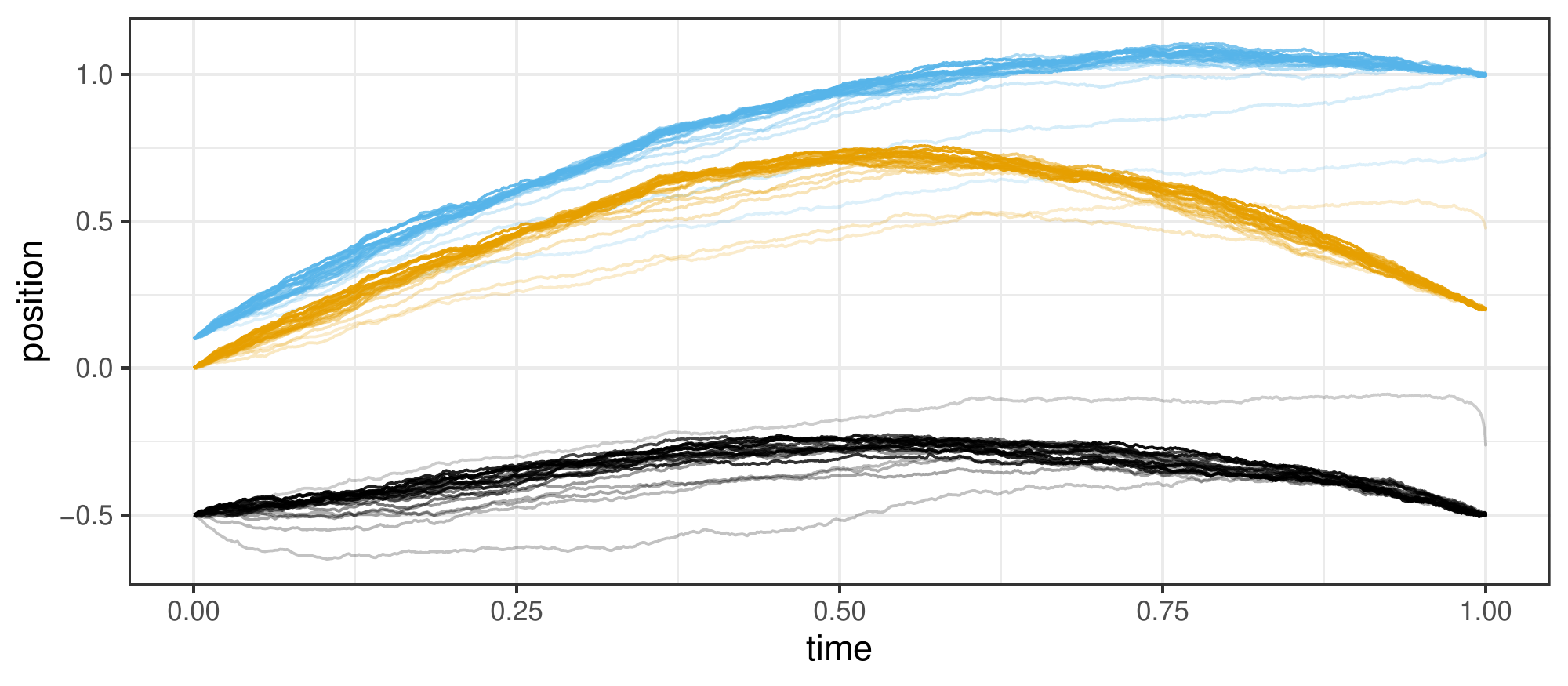}
\caption{Bridges for the 1D landmarks example, Eulerian model. The sampler ran for  $20\_000$ iterations. Every $1000$-th iterations is shown. ``Early'' iterations are depicted less thick.  Note that the initial trajectory does not satisfy the conditioning.\label{fig:bridges_1D_ahs}}
\end{center}
\end{figure}

\begin{figure}
\begin{center}
\includegraphics[scale=0.6]{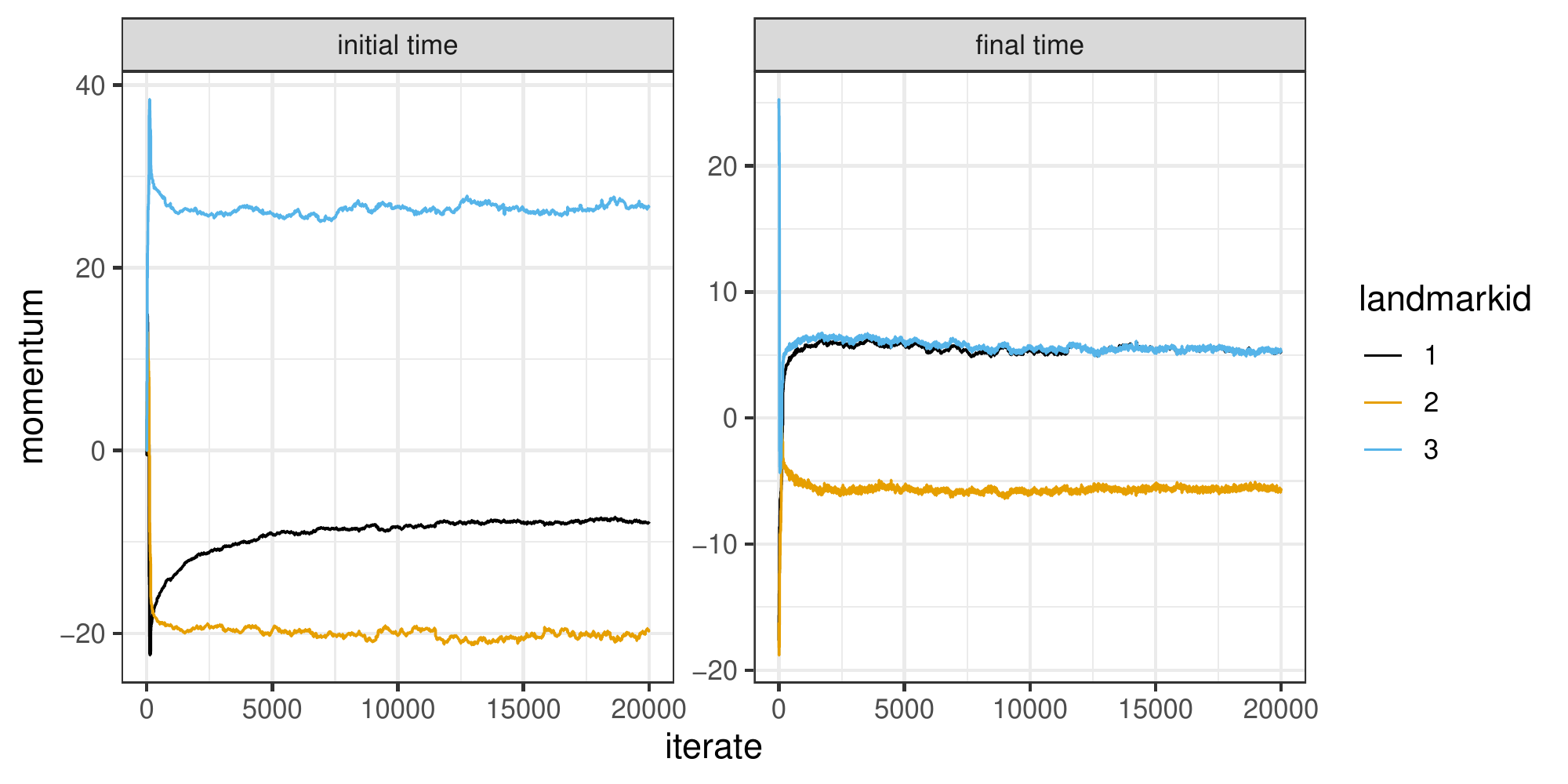}
\caption{Lagrangian model, 1D landmarks example. Traceplots showing iterates of momenta at times $0$ and $1$.Every $10$-th iterations is shown.\label{fig:momenta_1D_ms}}
\end{center}
\end{figure}

\begin{figure}
\begin{center}
\includegraphics[scale=0.6]{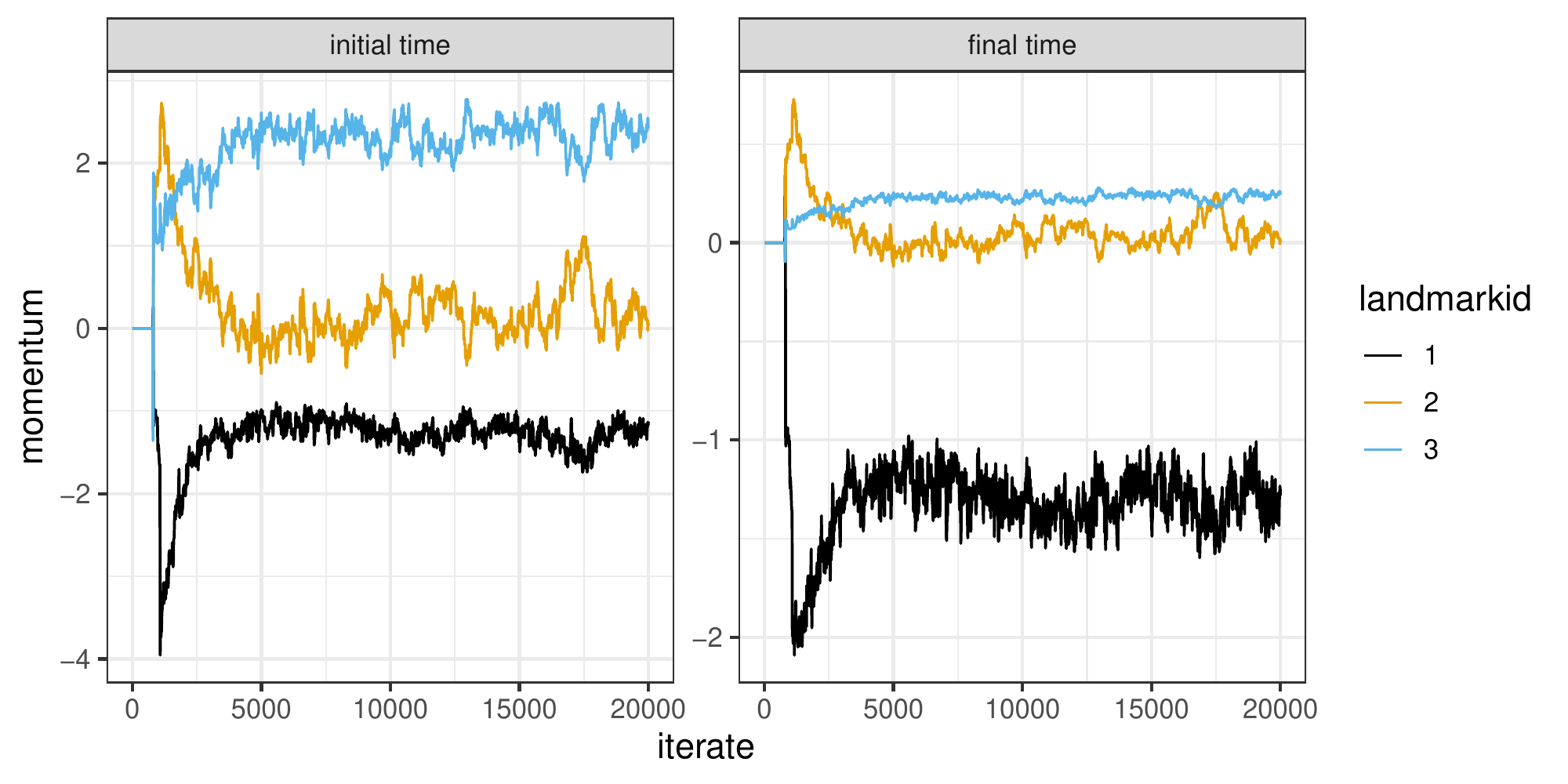}
\caption{Eulerian model, 1D landmarks example. Traceplots showing iterates of momenta at times $0$ and $1$.Every $10$-th iterations is shown.\label{fig:momenta_1D_ahs}}
\end{center}
\end{figure}

Sampled bridges for the Euleria and Eulerian model are in Figure \ref{fig:bridges_1D_ms} and Figure \ref{fig:bridges_1D_ahs} respectively. Clearly, upon initialisation the bridges do not match the conditionings at time $t=1$. This  is however resolved within a couple of thousand iterations. In Figures \ref{fig:momenta_1D_ms} and \ref{fig:momenta_1D_ahs} we show for both models traceplots of the momenta, both at the initial time and final time (note that we have matched the colours in the four figures). This foremost serves as a way to detect whether or not the chain has reached its stationary regime. Secondly, it is interesting to observe that for the Lagrangian model, the initial momenta are quite high and that for the upper two curves (blue and orange) the initial momenta are aposteriori in opposite directions. This can be explained by the value of the Hamiltonian kernel parameter (taken equal to $1$), which is deliberately chosen relatively large, forcing the orange and blue paths te remain close. However, as we condition the paths at distant locations ($0.2$ and $1.0$), the opposite initial momenta ensure the conditionings at time $1$ to be satisfied. 



\subsubsection{2D example from \cite{joshi_landmark_2000}}
Here we consider an example with $12$ $2$-dimensional landmarks from \cite{joshi_landmark_2000}. In this case we also include the Hamiltonian kernel parameter $a$ in the sampler, which is endowed with the Pareto distribution, More precisely, we assume $p(a) = 0.1 a^{-2} \ind_{[0.1,\infty)}(a)$ (hence the shape and scale parameters are taken $1$ and $0.1$ respectively). Having $a$ bounded away from zero is beneficial for numerical stability. 
Just like the previous example, for the Lagrangian model we fixed  $\gamma=1/\sqrt{n}\approx 0.29$, for the Eulerian model we fixed $\gamma=0.1$. 

\begin{figure}
\begin{center}
	\includegraphics[scale=0.6]{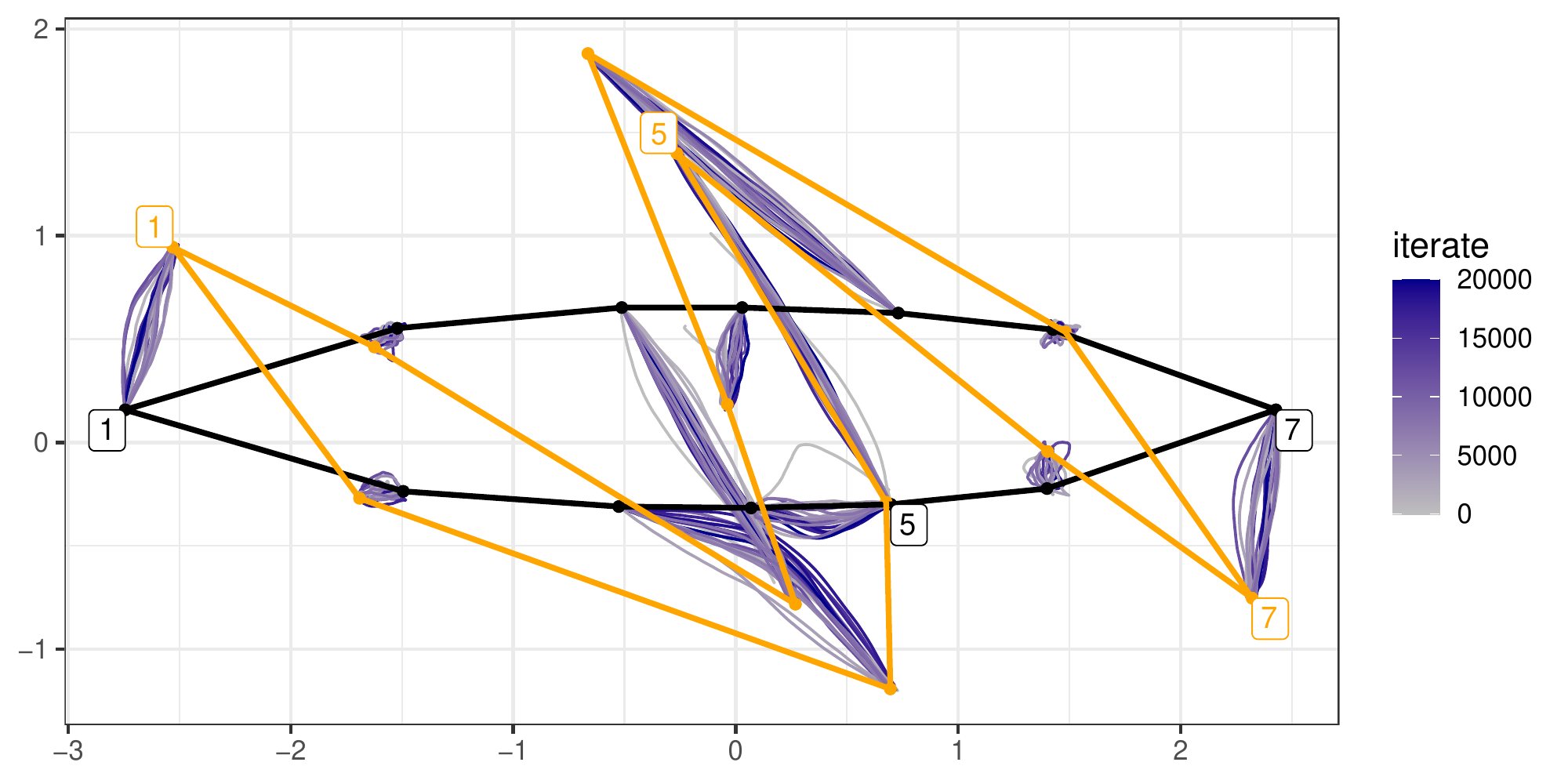}
\caption{Bridges for the Joshi-Miller example, Lagrangian model. The sampler ran for  $20\_000$ iterations. Every $1000$-th iterations is shown.  \label{fig:bridges_millerjoshi_ms}}
\end{center}
\end{figure}

\begin{figure}
\begin{center}
	\includegraphics[scale=0.6]{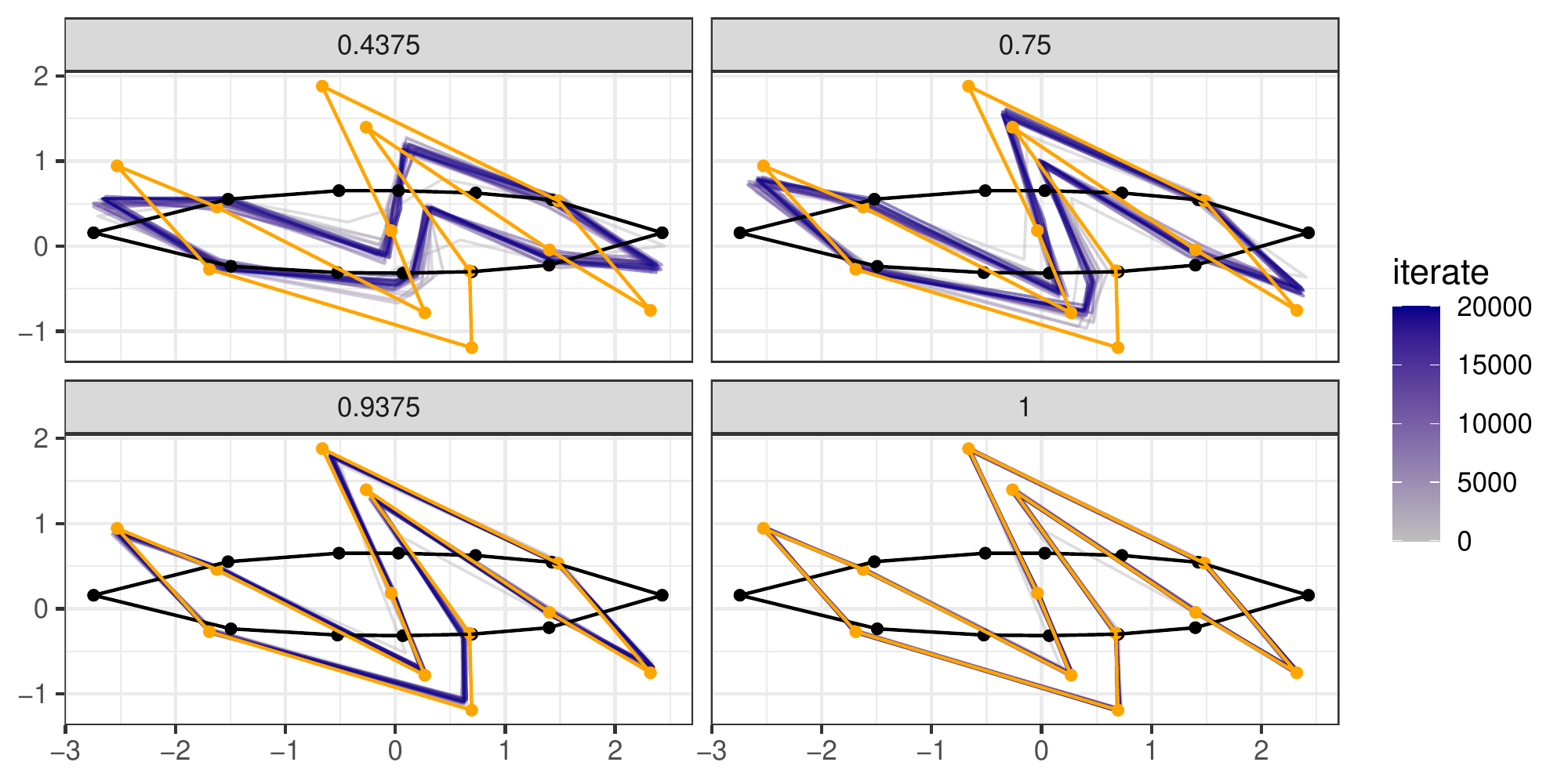}
\caption{Shape evolution  for the Joshi-Miller example, Lagrangian model. The sampler ran for  $20\_000$ iterations. Every $1000$-th iterations is shown. \label{fig:shapes_evolution_millerjoshi_ms}}
\end{center}
\end{figure}

\begin{figure}
\begin{center}
\includegraphics[scale=0.6]{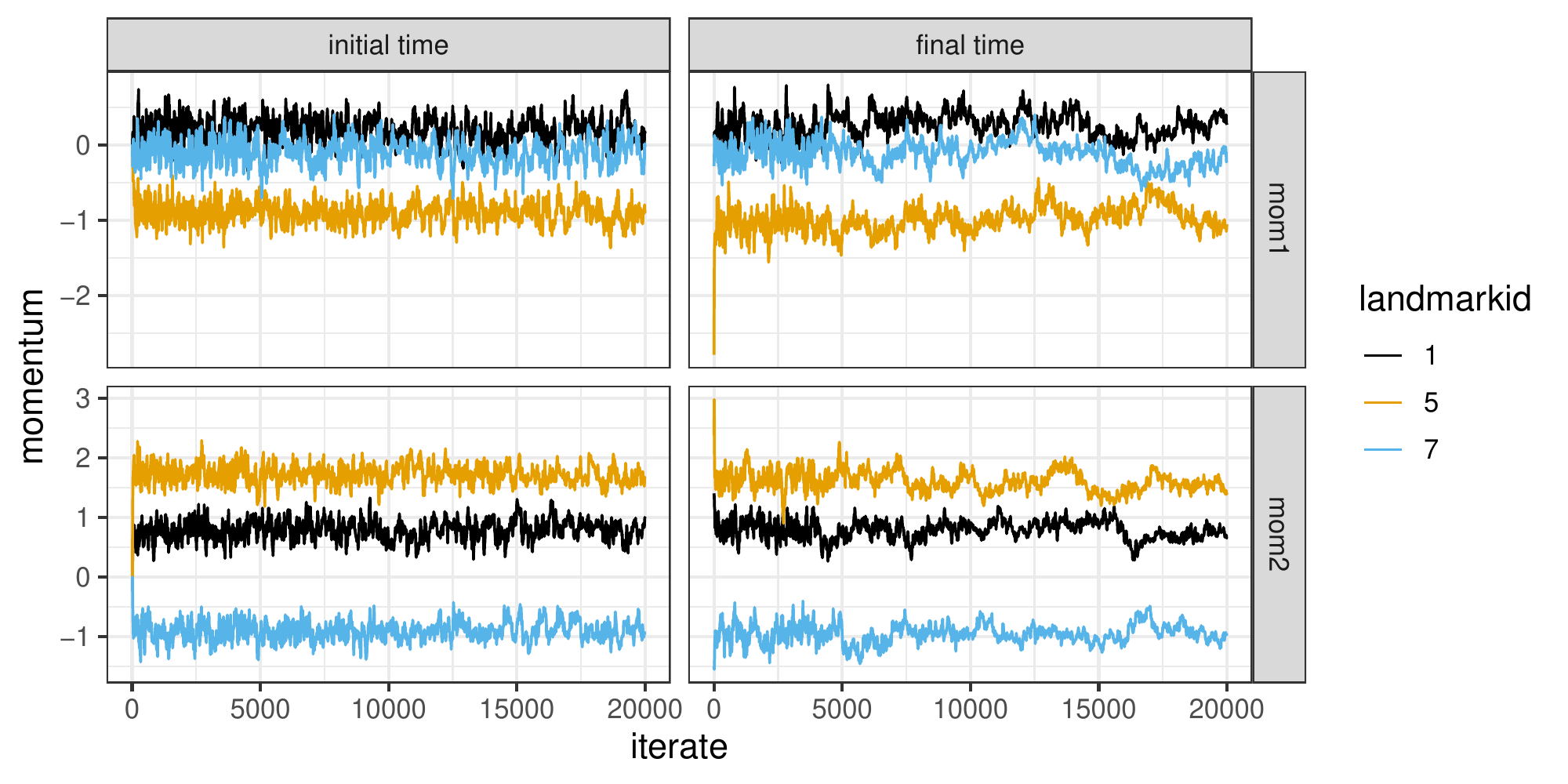}
\caption{Joshi-Miller example, Lagrangian model. Traceplots showing iterates of momenta at times $0$ and $1$. Every $10$-th iterations is shown.\label{fig:momenta_millerjoshi_ms}}
\end{center}
\end{figure}

\begin{figure}
\begin{center}
\includegraphics[scale=0.6]{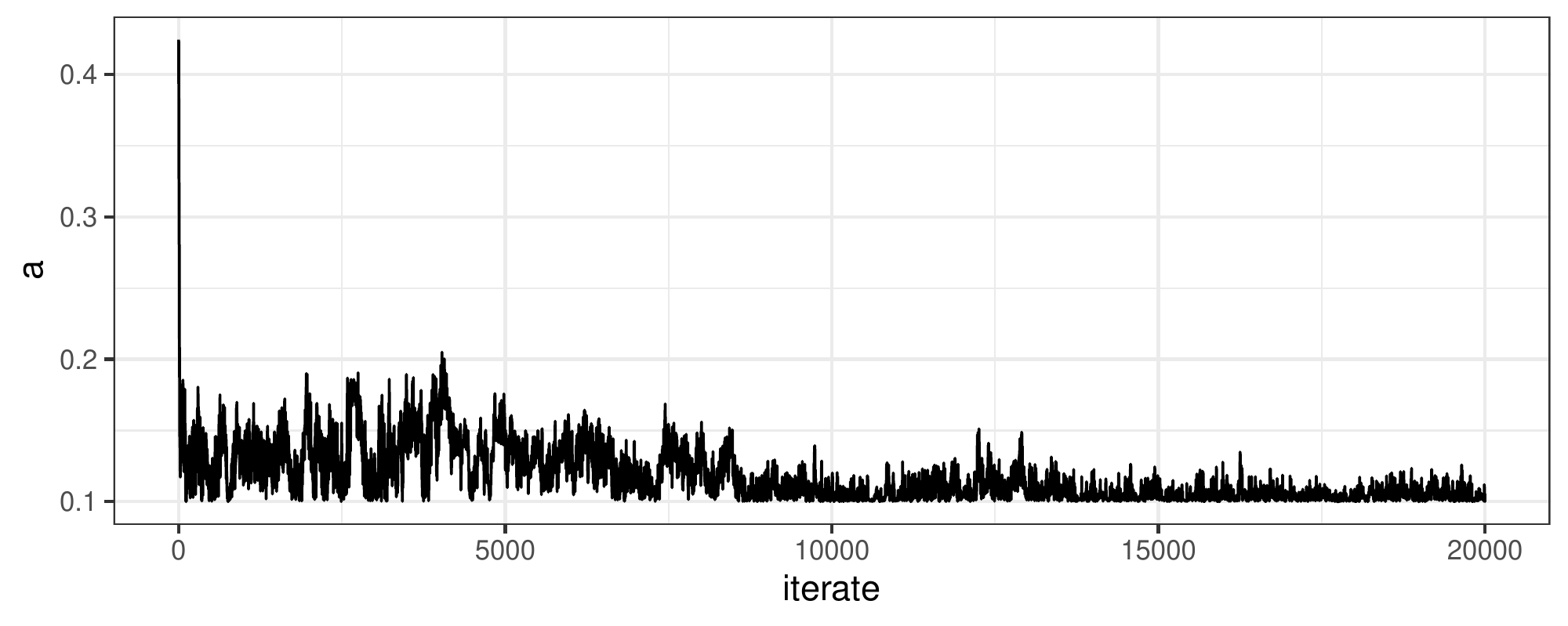}
\caption{Joshi-Miller example, Lagrangian model. Traceplots showing iterates of the Hamiltonian kernel parameter $a$. \label{fig:trace_millerjoshi_a_ms}}
\end{center}
\end{figure}

We ran a MCMC-sampler with $20\_000$ iterations, where one iterate consists of applying one step of Algorithms \ref{alg1}, \ref{alg2} and \ref{alg3}. 
We include figures for the Lagrangian model. In Figure \ref{fig:bridges_millerjoshi_ms} and in Figure \ref{fig:shapes_evolution_millerjoshi_ms} we visualise bridges and shape-evolution. Figure \ref{fig:momenta_millerjoshi_ms} shows traceplots for both coordinates of the momentum at  the initial time and final time. From these figures one can deduce for example that for landmark number $5$ the initial momentum is quickly increased after the initial iterations, to ensure that the bridge path indeed moves in north-direction. Figure \ref{fig:trace_millerjoshi_a_ms} gives a traceplot of the Hamiltonian kernel parameter $a$.  One can see that it takes about $10\_000$ iterations for the sampler to stabilise, and the corresponding change in sampled bridges is quite apparent in Figure  \ref{fig:bridges_millerjoshi_ms} (switch from greyish to purple-coloured paths).



\subsubsection{Corpus callosum} 
\begin{figure}
	\begin{center}
		\includegraphics[scale=0.6]{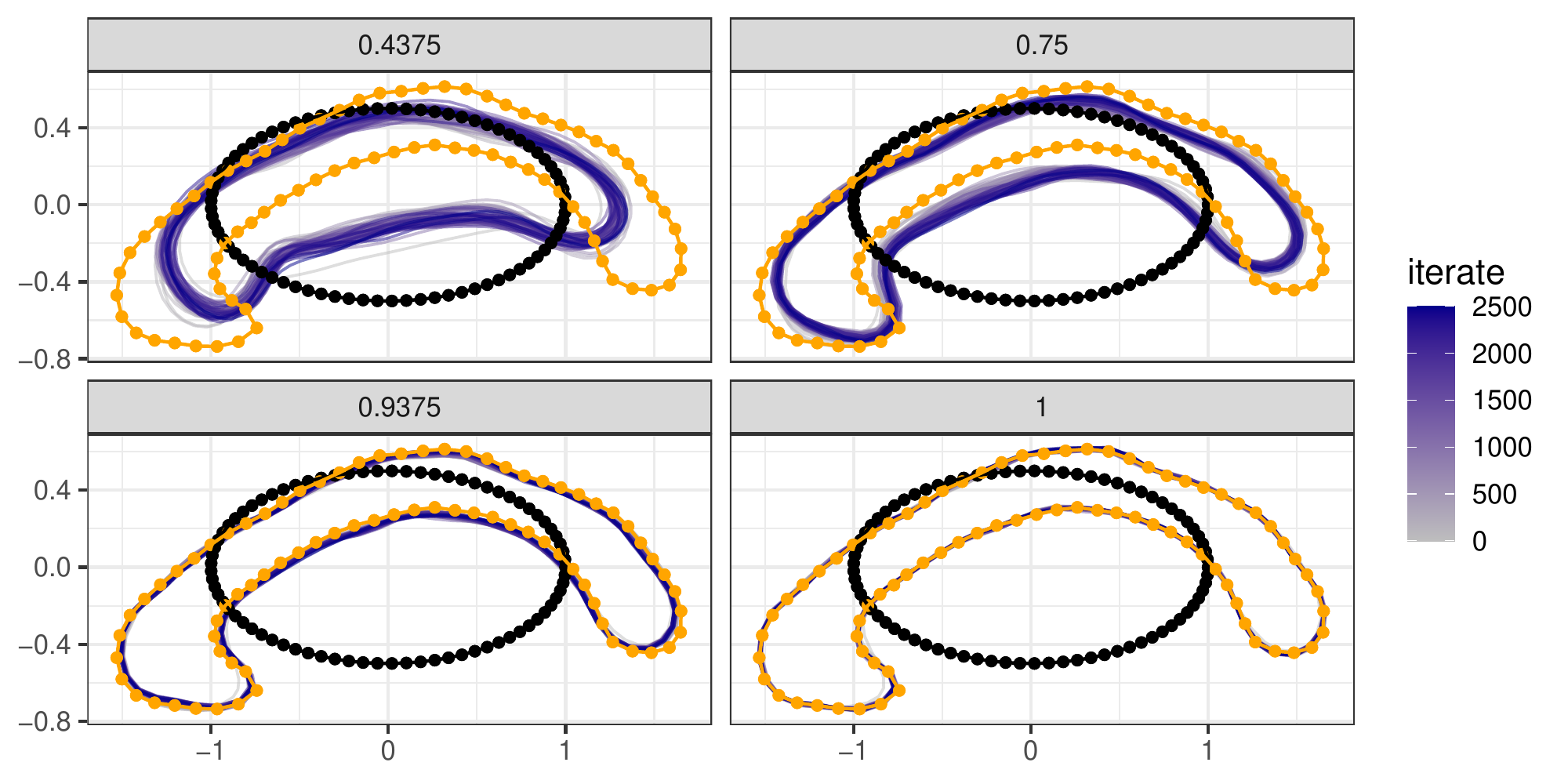}
	\end{center}
	\caption{Corpus callosum data. Initial shape (black) and final shape (orange) with $77$ landmarks. Every $10$th iteration from $500$ iterations is shown.\label{fig:shapes_evolution_dgs_ms}}
\end{figure}

\begin{figure}
	\begin{center}
		\includegraphics[scale=0.6]{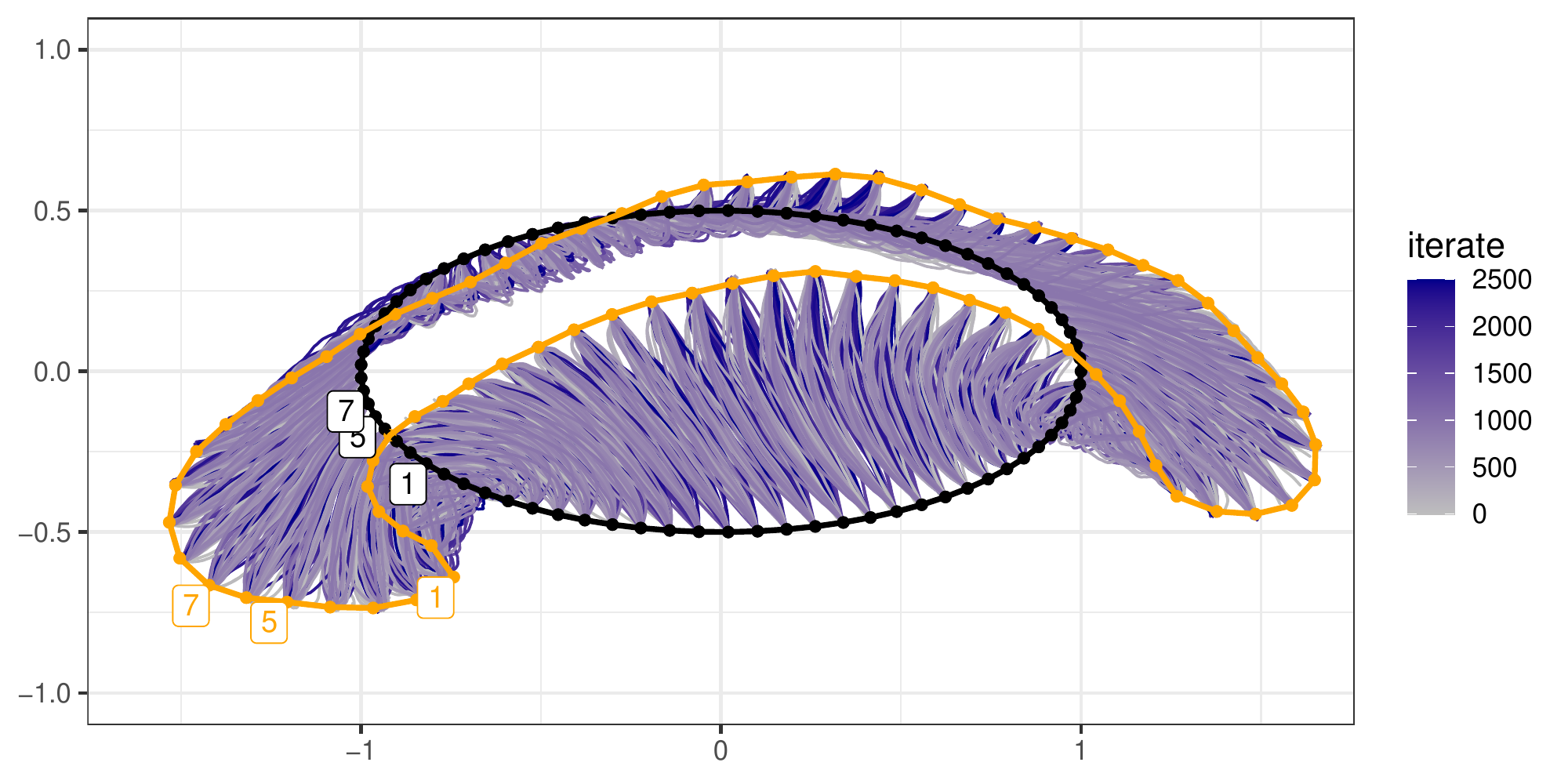}
	\end{center}
	\caption{Shape evolution  for the Corpus callosum data. The sampler ran for  $500$ iterations. Every $10$-th iterations is shown. \label{fig:bridges_dgs_ms}}
\end{figure}

\begin{figure}
	\begin{center}
		\includegraphics[scale=0.6]{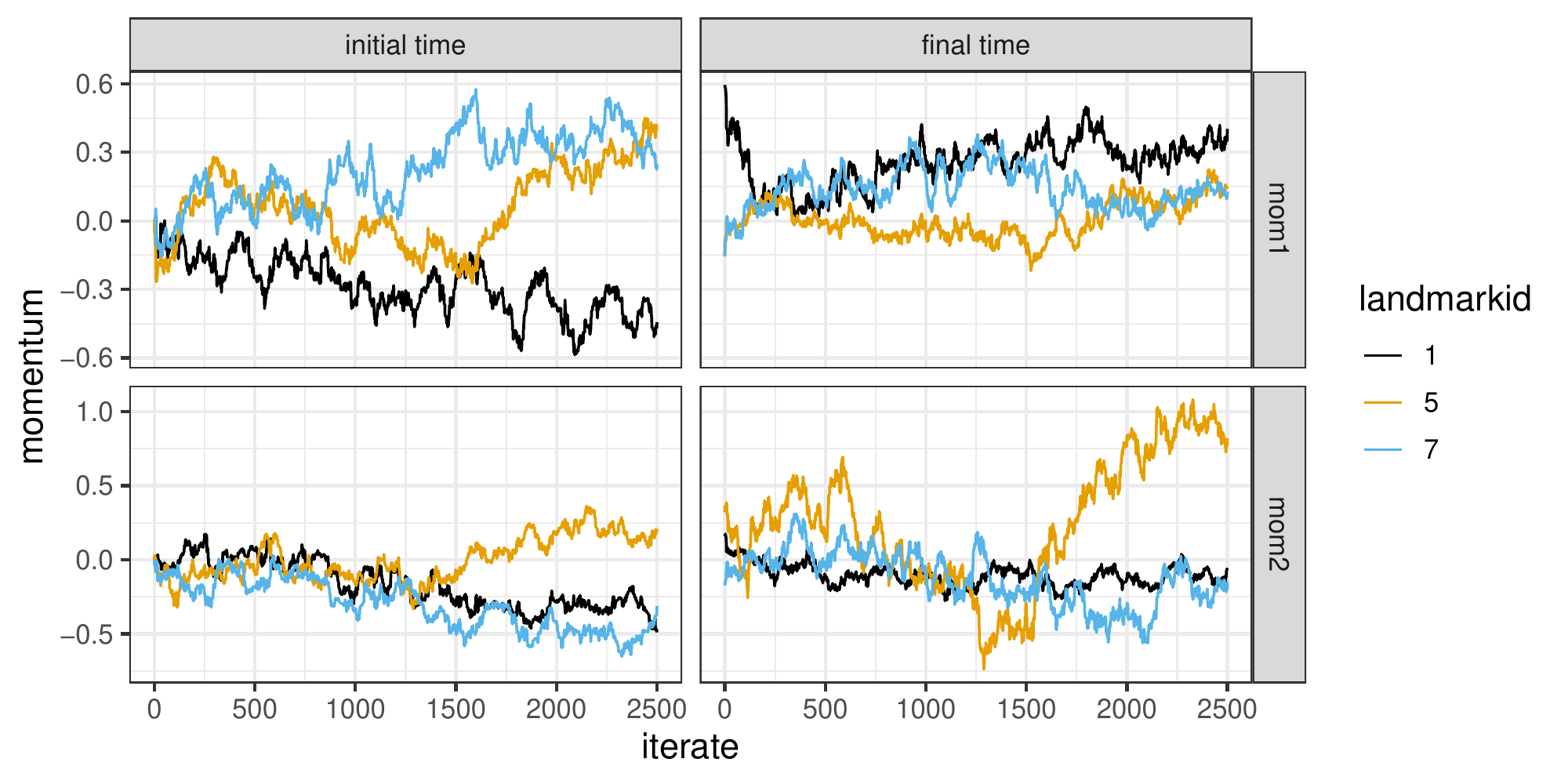}
	\end{center}	
	\caption{Corpus callosum data.Traceplots showing iterates of momenta at times $0$ and $1$. \label{fig:momenta_dgs_ms}}
\end{figure}

Here, we consider an example with many more landmarks where we start from an ellipse and transform to landmarks representing points on the outline of a human corpus callosum. In total there are $77$ landmarks. Here we used the Lagrangian model with 
 $a=0.2$ (Hamiltonian kernel parameter)  the average distance between adjacent landmarks on the shape and $\gamma)=1/\sqrt{n}\approx 0.12$. 
 In Figure \ref{fig:bridges_dgs_ms} we show sampled bridges over $2500$ iterations and in Figure  \ref{fig:shapes_evolution_dgs_ms} correspondingly the deformed shapes at 4 time instances in $[0,1]$. Traceplots of the momenta at times $0$ and $1$ are shown in Figure \ref{fig:momenta_dgs_ms}. 
 This example illustrates that the proposed algorithms can also be used in settings with a large number of landmarks. Note that after about $200$ iterations the bridges visually look very similar, whereas the early (light-greyish) paths look different.


\subsection{Examples for template estimation}

\subsubsection{Simulated shapes on an ellipse}\label{sec:template_ellipse} 

Here, we simulated forward $10$ trajectories   $n=15$ landmarks where the initial configuration is an ellipse. 
In the MCMC estimation, we took the Lagrangian model  and fixed the initial momenta to be equal to zero but include parameter estimation. Hence, each MCMC-iteration consists of executing one step of Algorithms \ref{alg1}, \ref{alg3} and \ref{alg:rmmala}. The parameter $a$ is endowed with the Pareto distribution, More precisely, we assume $p(a) = 0.1 a^{-2} \ind_{[0.1,\infty)}(a)$ (hence the shape and scale parameters are taken $1$ and $0.1$ respectively). 
We  ran the sampler for $2500$ iterations, saving every $50$th configuration.  We deliberately initialised the template configuration incorrectly by taking one of the observed configuration and rotating and stretching the shape. This is done to illustrate that the Riemannian manifold MALA steps perform satisfactory. In any practical application one could initialise the template shape by one of the observed shapes. 

\begin{figure}
	\begin{center}
		\includegraphics[scale=0.6]{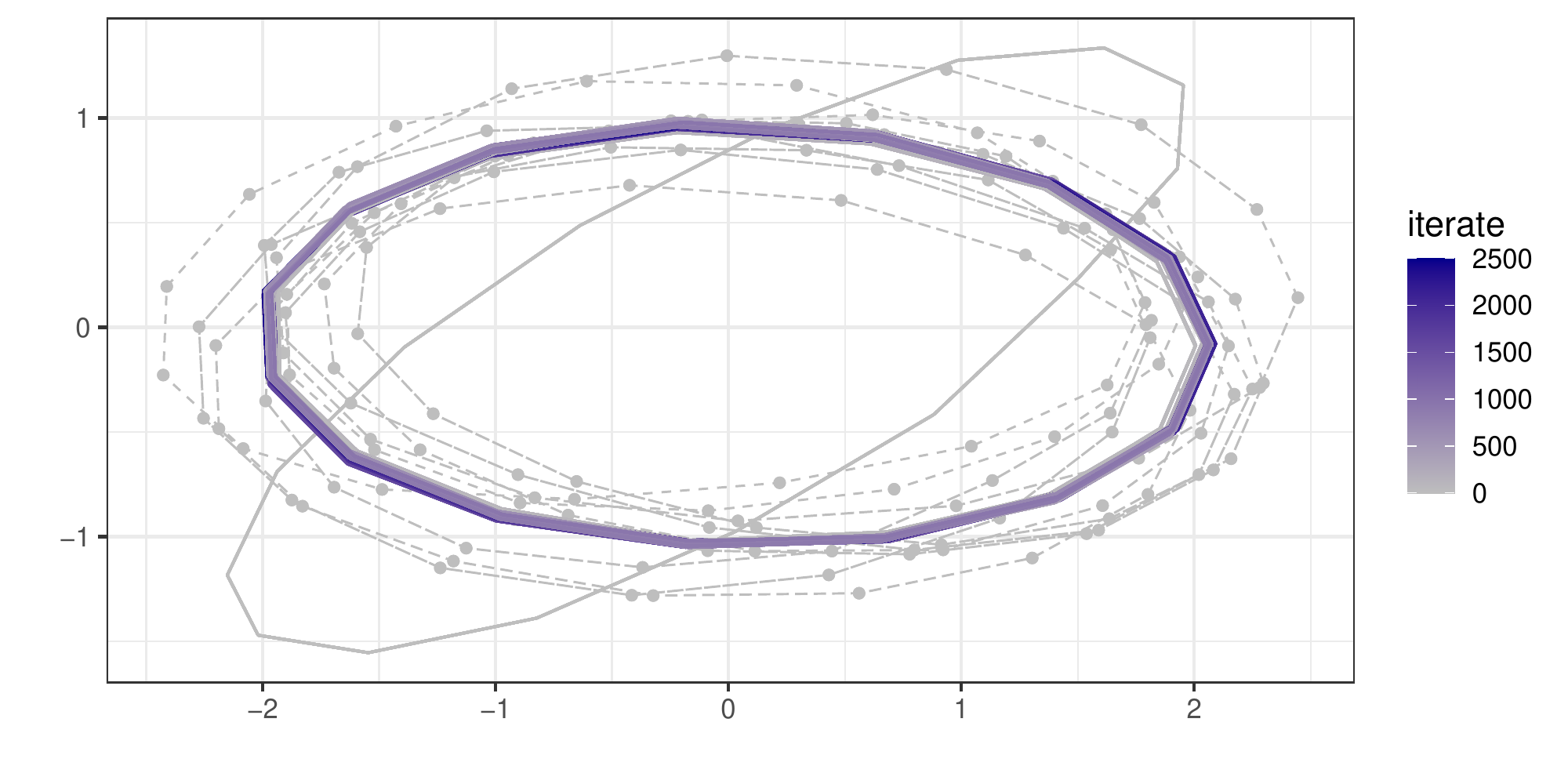}
	\end{center}
\caption{Initial shape configuration for the example in Section \ref{sec:template_ellipse} using the Lagrangian model.  The dashed grey curves are the $10$ observed configurations. Every $20$th MCMC-iterate is shown.  Note that that RMMALA quickly ``discovers'' the deliberately wrong initialisation. \label{fig:initial_shapes_exp2_ms}}
\end{figure}

\begin{figure}
	\begin{center}
		\includegraphics[scale=0.6]{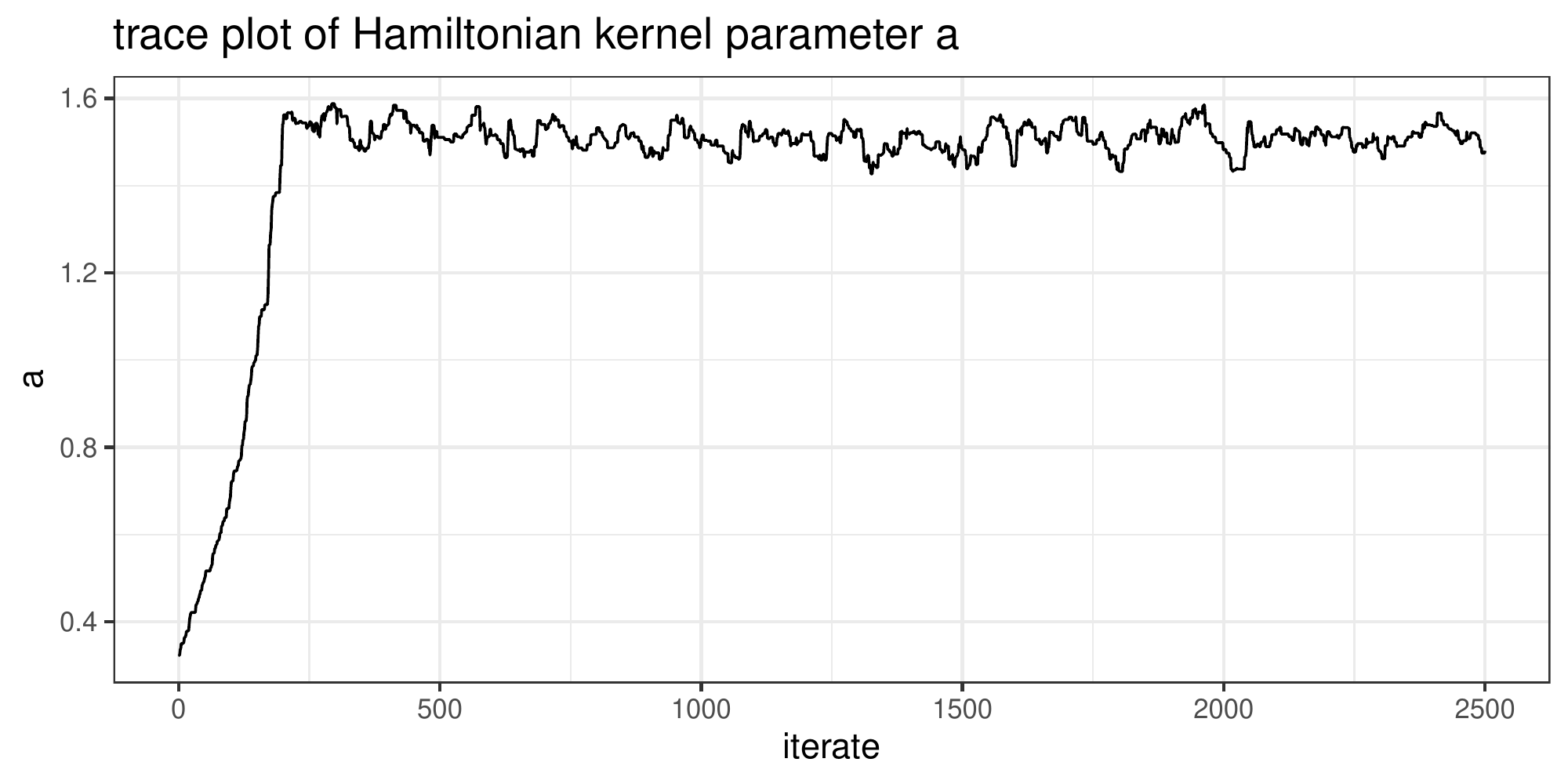}
	\end{center}
\caption{Trace plot for Hamiltonian kernel parameter for the example in Section \ref{sec:template_ellipse} using the Lagrangian model.   \label{fig:trace_exp2_a_ms}}
\end{figure}

Iterators of the template and parameter $a$ are shown in Figures \ref{fig:initial_shapes_exp2_ms} and \ref{fig:trace_exp2_a_ms} respectively. Clearly, the initial configuration can be recovered quite well, despite there are only $10$ observed shapes. From the traceplots of the parameter $a$ it appears there is  mild posterior uncertainty about its value. 

\subsubsection{Cardiac data}\label{subsec:cardiac}

\begin{figure}
	\begin{center}
		\includegraphics[scale=0.6]{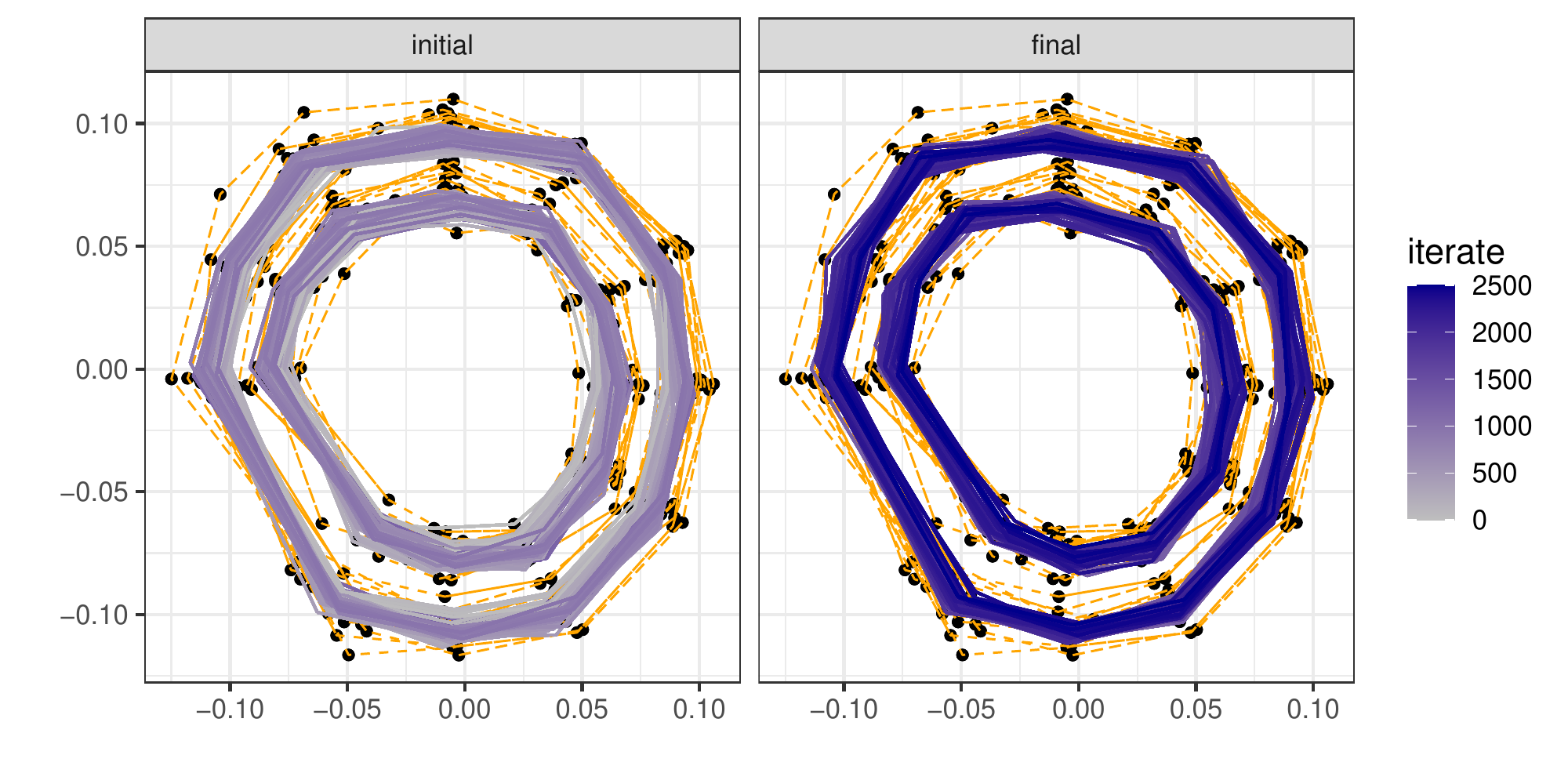}
		\includegraphics[scale=0.6]{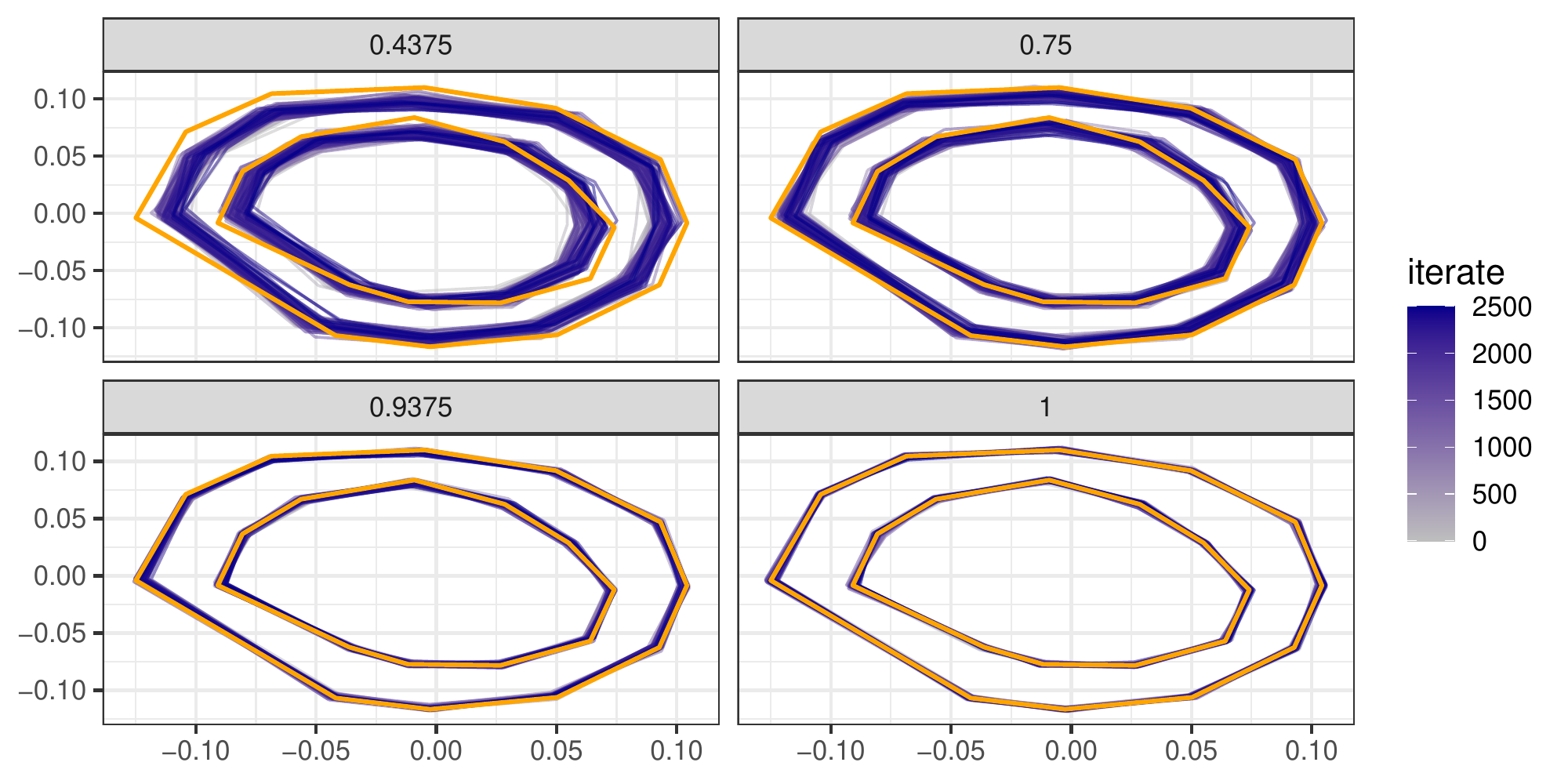}
				\includegraphics[scale=0.6]{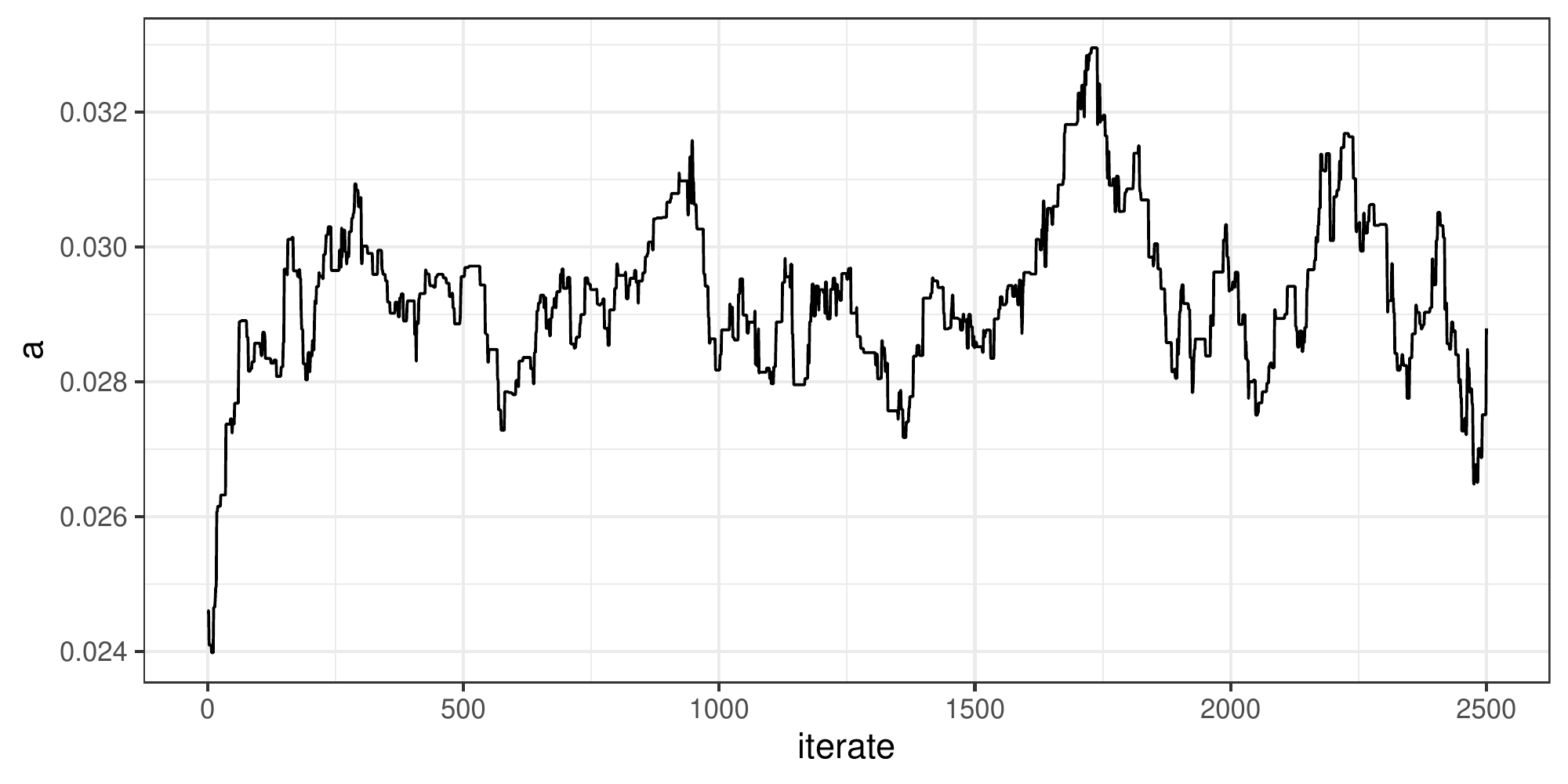}		
	\end{center}
	\caption{Cardiac data example from Section \ref{subsec:cardiac}. Top: iterates for template shape (locations of landmarks at time zero). The first 1250 iterates are in the left panel ``initial''; the final 1250 iterates in the right panel ``final''. Middle: Evolution over time of one chosen (fixed) shape (the iterates are shown at times $0.4375$, $0.75$, $0.9375$ and $1.0$). In both the top and middle figures the observed landmark locations are depicted by black dots, connected by orange line segments.  Bottom: trace plot for the Hamiltonian kernel parameter.\label{fig:cardiac_initialshape}}
\end{figure}

 In Figure \ref{fig:cardiac_initialshape} we show the results when  applying the algorithm to landmarks on 14 cardiac images of human left ventricles~\cite{stegmann_extending_2001}. The manually annotated landmarks are consistently placed over the set of shapes. Based on $22$ landmarks for each shape we estimate the template shape and Hamiltonian kernel parameter. This experiment extends the template estimation with Brownian motion performed in~\cite{sommer_bridge_2017} to the models in focus in this paper. We used parameter updating (assuming apriori $p(a) = 0.01 a^{-2} \ind_{[0.01,\infty)}(a)$), bridge updating and initial state updating (using RMMALA), while fixing the initial momenta to zero. Hence, just as the previous subsection, each MCMC-iteration consists of executing one step of Algorithms \ref{alg1}, \ref{alg3} and \ref{alg:rmmala}. We used $2500$ iterations. From Figure \ref{fig:cardiac_initialshape} it appears that the chain reaches its stationary region within $500$  iterations. 

\section{Discussion and future work}

We have extended the framework of~\cite{schauer_guided_2017,partialbridges} for sampling diffusion bridges to nonlinear high dimensional stochastic Hamiltonian systems. For this work, problems in computational anatomy~\cite{younes_shapes_2010} revolving around systems of particles to describe images have been our main motivation. We used these models throughout to  illustrate the  application of diffusion bridges specifically in stochastic landmark dynamics. 
We left open many possible improvements of this method  such as allowing for multiple observations over time, resolving unknown landmark correspondence, a step towards infinite dimensional shape matching and the problem of the parametrisation of these shapes~\cite{bauer2014overview}. 

This work opens the door to not not only a wider use in statistical inference in computational anatomy, but other related fields, where high dimensional stochastic Hamiltonian systems are involved in the modelling, such as in biology, data assimilation, swarming, etc.
From a more theoretical point of view, various extensions of this work could be considered. The driving  Brownian motion in the SDE could for example be generalised  to a L\'evy process or even rough paths.

From a numerical perspective, we used classical methods for simplicity, but the presented algorithms could be extended to include more modern numerical integration schemes or MCMC methods, such as geometric integrators~\cite{milstein2002,hairer2006geometric,Niu2018,holm2018stochastic,buckwar2019spectral}, to preserves the geometrical structure of the problem, or some more advanced geometrical Hamiltonian Monte-Carlo methods~\cite{barp2018geometry,barp2019hamiltonian}. 
In addition, the use of automatic differentiation, similar as in \cite{kuhnel2017computational,kuhnel2019differential}, has turned out to be beneficial for the  implementation of efficient methods, and surely deserves to be more exploited for future works.

\section*{Acknowledgement} 
A substantial part of this research was carried out during a sabbatical leave of the second author at the University of Copenhagen (UCPH). Financial support by the faculty of Science, UCPH trough the programme ``Visiting scholars at SCIENCE'' and by both the mathematics and computer science department of UCPH are greatly acknowledged. The work presented is in addition supported by the CSGB Centre for Stochastic Geometry and Advanced Bioimaging funded by a grant from the Villum foundation, the Villum Foundation grant 00022924, and the Novo Nordisk Foundation grant NNF18OC0052000. AA acknowledge funding through EPSRC award EP/N014529/1 supporting the EPSRCCentre for Mathematics of Precision Healthcare at Imperial.

\appendix

\section{Markov Chain Monte Carlo}\label{sec:mh}

A \emph{Markov chain Monte Carlo (MCMC)} method for sampling from  a distribution $\pi$ is any method producing an ergodic Markov chain whose stationary distribution is $\pi$. Virtually all such algorithms can be seen of specific instances of the Metropolis-Hastings algorithm. For this paper it suffices to assume that $\pi$ has a density with respect to Lebesgue measure on $\RR^d$. We will denote this density by $\bar\pi$. Suppose $Q$ is a Markov kernel on $\RR^d$ with density $q$. That is, for measurable sets $B$, $x\mapsto Q(x,B)$ is a measurable mapping and for fixed $x$, $B\mapsto Q(x,B)=\int_B q(x,y) \dd y$ is a probability measure.
\begin{defn}\label{def:mhalg}
The \emph{Metropolis-Hastings (MH) algorithm} is the algorithm by which a Markov chain is constructed which evolves  $x_n=x$ to $x_{n+1}$ by the following steps
\begin{enumerate}
\item propose $y$ from the proposal kernel $Q(x,\cdot)$;
\item 
Compute $$ \alpha(x,y)=\min \left(1, \frac{\bar\pi(y)}{\bar\pi(x)} \frac{q(y,x)}{q(x,y)}\right). $$
\item
Set 
\[ x_{n+1} = \begin{cases} y & \quad \text{with probability $\alpha(x,y)$} \\ 
 x & \quad \text{with probability $1-\alpha(x,y)$} \end{cases} . \]
\end{enumerate}
\end{defn}
Under weak assumptions, this algorithm indeed produces an ergodic Markov chain with stationary distribution  $\pi$. 
For computing $\alpha(x,y)$ it suffices  to know $\pi$ up to a proportionality constant. Within Bayesian statistics, this is a very attractive property of the algorithm, as the posterior is typically only known up to an integration constant. 

The specification of the algorithm leaves open the choice of the proposal kernel $Q$. A symmetric random walk proposal $q(x,y) = \phi(y-x)$, with $\phi$ a symmetric density around zero, is a simple choice that may work if the dimension $d$ is small. Otherwise, a more sophisticated proposal kernel that takes gradient information into account may be computationally more efficient. The Metropolis Adjusted Langevin Algorithm (MALA) is based on the fact that the solution to the stochastic differential equation
\[ \dd X_t = A \nabla  \log \bar\pi(X_t) \dd t + \sqrt{2 A} \dd W_t, \]
with $W$ denoting Brownian Motion, has invariant distribution $\pi$. Here, $A$ is a strictly positive definite preconditioning matrix and  $\sqrt{A}$ is the matrix-square root of $A$. 
 Upon discretising this equation with stepsize $\tau$ (Euler forward discretisation) we obtain the proposal density $q(x,y) = \psi(y; x+  A \nabla  \log \bar\pi(x), 2A)$ with (as before) $\psi(x; \mu,\Sigma)$ denoting the density of the $N(\mu, \Sigma)$-distribution, evaluated at $x$. Both $\tau$ and $A$ are tuning parameters; $\tau$ acting as a stepsize. In general, the choice of $A$ is dictated by the geometrical properties of the density $\bar\pi$. This is the basic idea of the Riemannian MALA algorithm. 
If the dimension $d$ is large, one often uses a proposal where only  parts for the current iterate $x_n=x$ are stochastically perturbed. This is sometimes called ``one-at-a-time'' sampling, ``substitution sampling'', or simply ``Gibbs sampling''. Clearly, this is just a special case of the Metropolis-Hastings algorithm. 

\medskip

For illustration purposes, we give a  finite-dimensional version corresponding to Algorithm \ref{alg1} (updating of bridges for fixed parameter and initial state). Hence assume a data-generating model of the form $w \sim N(0, I)$, $x = F(w)$ leading to the observation  $y = h(x) + \epsilon$, $\epsilon$ being a random quantity independent of $w$. Assume the likelihood for the parameter $x$ is given by  $L(x)=L(F(w))$. Suppose $L(F(w)) = c \Psi(G(w))$, for fixed maps $G$ and $\Psi$, the constant being independent of $w$. Using the proposal density $q(w^\circ \mid w)$,  the acceptance probability for the Metropolis-Hastings algorithm is given by $A \wedge 1$, where $ A = \frac{q(w \mid w^\circ)}{q(w^\circ \mid w)}   \frac{\psi(w^\circ; 0, I)}{\psi(w; 0, I)}  \frac{\Psi(G(w^\circ))}{\Psi(G(w))} $.  
If $q$ is chosen according to the Crank-Nicolson rule, i.e.\ $w^\circ \mid w \sim N(\eta w, (1-\eta^2)I)$, then the first two terms in the expression for $A$ cancel and we end up with 
$ A =  \Psi(x^\circ)/\Psi(x)$, 
where $x^\circ = G(w^\circ)$ and $x = G(w)$. 
This expression remains valid if $w$ is a Wiener-process and $G=\GP$ is a guided proposal. 



\section{Matching conditions on the auxiliary process for exact landmarks matching}\label{app:matching_conditions}
 For exact matching, absolute continuity will only hold if the parameters of the  auxiliary process  $\widetilde{X}$,  satisfy certain {\it matching conditions}.    
In case the diffusion is uniformly elliptic,  $L_T=I_{N\times N}$,  it was shown in~\cite{schauer_guided_2017} that $\widetilde{a}=\si\sigma^\T$ must satisfy $\widetilde{a}(T)=a(T,x_T)$  for absolute continuity.   In the more general case considered here, the matching conditions look somewhat more difficult, a precise statement is given in section 2 of~\cite{partialbridges}. Especially for the Eulerian model, these conditions are hard to verify (mainly due to the fact that the diffusivity is state dependent and not all components of  $X_T$ are observed). 
However, the results in this paper suggest that  $\PP^\star$ is absolutely continuous with respect to $\PP^\circ$  if {\it (i)}  condition  \eqref{eq:match-conject} is satisfied; {\it (ii)} the number of noise-sources $J$ satisfies  $J\ge Nd$ ($N$ and $d$ denoting the number of landmarks and their dimension respectively).
 This conjecture is numerically confirmed by experiments in Section 4.1 of~\cite{partialbridges}. The second condition that requires sufficiently many noise sources has also been discussed in Remark 4.1 in  \cite{arnaudon_geometric_2019}.

\section{Proof of Proposition \ref{lem:ahs-strat2ito} on Stratonovich to It\^{o} correction for Eulerian model}

In coordinates, the stochastic equations for the Eulerian model are
\begin{align}
    \begin{split}
      d q_i^\alpha &= \frac{\partial h}{\partial p_i^\alpha} dt + \sum_{l=1}^J \sigma_l^\alpha(q_i) \circ d W_t^l   \,, \\
    d p_i^\alpha &= -\frac{\partial h}{\partial q_i^\alpha} dt - 
    \sum_{l=1}^J\sum_\beta\frac{\partial\sigma_l^\beta(q_i)}{\partial q_i^\alpha } p_i^\beta \circ  d W_t^l \, ,
    \end{split}
    \label{sto-Ham-coor}
\end{align}
The same process with It\^o's integrals, has the additional term 
\begin{equation}\label{eq:extraq}
        dq_i^\alpha += \frac12\frac{\partial\sigma_l^\alpha(q_i) }{\partial q_i^\beta} \sigma_l^\beta(q_i) \, ,
\end{equation}
in the $q$ equation and for the $p$ equation 
\begin{equation}\label{eq:extrap}
        dp_i^\alpha += \frac12 p_i^\gamma \frac{\partial\sigma_l^\gamma(q_i)}{\partial q_i^\beta}\frac{\partial\sigma_l^\beta(q_i)}{\partial q_i^\alpha} - \frac12 p_i^\beta \frac{\partial^2\sigma_l^\beta(q_i) }{\partial q_i^\alpha \partial q_i^\gamma} \sigma_l^\gamma(q_i)\, .
\end{equation}

The extra term for $q^\alpha$ given in \eqref{eq:extraq} equals 
\begin{align*} \frac12 \sum_\ell \sum_\beta \frac{\partial \sigma_\ell^\alpha(q)}{\partial q^\beta} \sigma_\ell^\beta(q) &=\frac12 \sum_\ell \sum_\beta \gamma_\alpha \nabla_\beta \bar{k}_\tau(q-\delta) \gamma_\beta \bar{k}_\tau(q-\delta)\\ & =  \frac12 \sum_\ell \langle \nabla \bar{k}_\tau(q-\delta), \gamma\rangle \bar{k}_\tau(q-\delta) \gamma_\alpha\, ,
\end{align*}
which gives \eqref{eq:stratito-q}. Notice that we write $\nabla_\beta$ to denote $(\partial)/(\partial q^\beta)$. 

For $p^\alpha$ the first term on the right-hand-side of equation \eqref{eq:extrap} is given by 
\begin{align*}
 \frac12 \sum_\ell \sum_u \sum_\beta p^u \frac{\partial\sigma_l^u(q)}{\partial q^\beta}\frac{\partial\sigma_l^\beta(q)}{\partial q^\alpha} &= \frac12 \sum_\ell \sum_u \sum_\beta  \gamma_u \gamma_\beta \nabla_\beta \bar{k}_\tau(q-\delta) \nabla_\alpha \bar{k}_\tau(q-\delta)  p^u \\&= \frac12 \sum_\ell \nabla_\alpha \bar{k}_\tau(q-\delta) \langle \gamma,p\rangle \langle \gamma, \nabla \bar{k}_\tau(q-\delta)\rangle\, .
\end{align*}
This gives the first term in \eqref{eq:stratito-p}.
The second term on the right-hand-side of equation \eqref{eq:extrap} is given by 
\begin{align*}
	- \frac12  \sum_\ell \sum_u \sum_\beta p^\beta \frac{\partial^2\sigma_l^\beta(q) }{\partial q^\alpha \partial q^u} \sigma_l^u(q) &= 
-\frac12 	 \sum_\ell \sum_u \sum_\beta p^\beta \gamma_\beta \left[ \nabla_\alpha\left( \nabla_u \bar{k}_\tau(q-\delta)\right) \right] \gamma_u \bar{k}(q-\delta) \\ & = -\frac12 	 \sum_\ell \langle(p,\gamma)\rangle \bar{k}_\tau(q-\delta)  \nabla_\alpha \left( \langle \nabla \bar{k}(q-\delta), \gamma \rangle \right)\, .
\end{align*}
This gives the second term in \eqref{eq:stratito-p}. 

\bibliographystyle{siamplain}
\bibliography{lit,ss}
\end{document}